\documentclass{article}
\usepackage[utf8]{inputenc}

\usepackage[titletoc]{appendix}
\usepackage{alltt}
\usepackage[T1]{fontenc}
\usepackage{titlesec, blindtext, color}
\usepackage{graphicx}
\usepackage{epstopdf}
\usepackage{amsfonts}
\usepackage{relsize}
\usepackage{soul}
\usepackage{marvosym}
\usepackage{amssymb,amsmath,amsthm}
\usepackage{dsfont, bbm}
\usepackage{mathtools}
\usepackage{mathrsfs}
\usepackage{footnote}
\usepackage{enumitem}
\usepackage{scalefnt}
\usepackage{xcolor}
\usepackage{adjustbox}
\usepackage{setspace}
\usepackage{tensor}

\usepackage{float}
\restylefloat{table}
\usepackage{array}
\usepackage{hvfloat}
\usepackage{fancyhdr}
\usepackage{wrapfig}
\usepackage{url}
\usepackage[nottoc,notlot,notlof]{tocbibind}

\usepackage{hyperref}
\hypersetup{
    pdftitle={Rough Fubini's theorem},
}
\usepackage[numbered]{bookmark}

\usepackage{tikz-cd}
\tikzcdset{every label/.append style = {font = \small}}

\allowdisplaybreaks
\titleformat{\section}
{\bfseries\large}{\thesection}{1em}{}

\usepackage[paperwidth=8.27in, paperheight=11.69in, top=1.2in, bottom=1.6in, outer=1.2in, inner=1.2in, heightrounded]{geometry}

\usepackage[backend=biber, style=alphabetic, maxbibnames=99, maxcitenames=50]{biblatex}

\newcommand{\Real}{\mathbb{R}}

\newcommand{\pfloor}{\lfloor p \rfloor}
\newcommand{\ptilfloor}{\lfloor \tilde{p} \rfloor}

\newcommand{\Xbf}{\mathbf{X}}
\newcommand{\Rbf}{\mathbf{R}}
\newcommand{\Vbf}{\mathbf{V}}

\newcommand{\Dcal}{\mathcal{D}}

\newcommand{\Ical}{\mathcal{I}}
\newcommand{\Pcal}{\mathcal{P}}

\newcommand{\Dscr}{\mathscr{D}}

\newcommand{\ptilde}{\tilde{p}}
\newcommand{\qtilde}{\tilde{q}}
\newcommand{\Atilde}{\tilde{A}}
\newcommand{\Itilde}{\tilde{I}}
\newcommand{\Ntilde}{\tilde{N}}
\newcommand{\Vtilde}{\tilde{V}}
\newcommand{\Xtilde}{\tilde{X}}
\newcommand{\Ytilde}{\tilde{Y}}
\newcommand{\Ztilde}{\tilde{Z}}
\newcommand{\Rtilde}{\tilde{R}}

\newcommand{\Dcaltilde}{\tilde{\mathcal{D}}}
\newcommand{\Rbftilde}{\tilde{\mathbf{R}}}

\newcommand{\omegatil}{\tilde{\omega}}
\newcommand{\thetatil}{\tilde{\theta}}

\DeclareMathOperator{\Hom}{\textnormal{Hom}}
\DeclareMathOperator{\Lip}{\textnormal{Lip}}
\DeclareMathOperator{\variation}{\textnormal{-var}}
\newtheorem{theorem}{Theorem}[section]

\newtheorem{condition}[theorem]{Condition}

\newtheorem{corollary}[theorem]{Corollary}

\newtheorem{definition}[theorem]{Definition}

\newtheorem{lemma}[theorem]{Lemma}

\newtheorem{remark}[theorem]{Remark}

\usepackage{subfiles}
\bibliography{RoughFubiniBib}

\title{\textbf{A Fubini type theorem for rough integration\thanks{The work of Thomas Cass is supported by EPSRC Programme Grant EP/S026347/1. The work of Jeffrey Pei is supported by the EPSRC Financial Computing and Business Analytics Centre for Doctoral Training}}}

\author{Thomas Cass\thanks{Imperial College London and the Alan Turing Institute (email: thomas.cass@imperial.ac.uk)} , Jeffrey Pei\thanks{Imperial College London and University College London (email: jeffrey.pei14@imperial.ac.uk)}}

\date{}

\begin{document}

\maketitle

\vspace{-3em}

\begin{abstract}\noindent
We develop the integration theory of two-parameter controlled paths $Y$ allowing us to define integrals of the form
\begin{equation}\label{joint integral}
    \int_{[s,t] \times [u,v]}
        Y_{r,r'}
    \;d(X_{r}, X_{r'})
\end{equation}
where $X$ is the geometric $p$-rough path that controls $Y$. This extends to arbitrary regularity the definition presented for $2\leq p<3$ in the recent paper of Hairer and Gerasimovi\v{c}s \cite{GH19} where it is used in the proof of a version of H\"{o}rmander's theorem for a class of SPDEs. We extend the Fubini type theorem of the same paper by showing that this two-parameter integral coincides with the two iterated one-parameter integrals 
\[
    \int_{[s,t] \times [u,v]}
        Y_{r,r'}
    \;d(X_{r}, X_{r'})
    =
    \int_{s}^{t}
        \int_{u}^{v}
            Y_{r,r'}
        \;dX_{r'}
    \;dX_{r}
    =
    \int_{u}^{v}
        \int_{s}^{t}
            Y_{r,r'}
        \;dX_{r}
    \;dX_{r'}.
\]
A priori these three integrals have distinct definitions, and so this parallels the classical Fubini's theorem for product measures. By extending the two-parameter Young-Towghi inequality of \cite{Tow02} in this context, we derive a maximal inequality for the discrete integrals approximating the integral in \eqref{joint integral}. As a further benefit we remove the technical assumption made in \cite{GH19} that $Y$ be smoothly approximable. We also extend the analysis to consider integrals of the form 
\begin{equation}\label{joint integral 2}
    \int_{[s,t] \times [u,v]}
        Y_{r,r'}
    \;
    d(X_{r}, \Xtilde_{r'})
\end{equation}
for possibly different rough paths $X$ and $\Xtilde$, and obtain the corresponding Fubini type theorem. We prove continuity estimates for these integrals in the appropriate rough path topologies. As an application we consider the signature kernel \cite{sigkernelsgoursat}, which has recently emerged as a useful tool in data science, as an example of a two-parameter controlled rough path which also solves a two-parameter rough integral equation.\\
\end{abstract}

\bookmark[page=1,level=1]{Abstract}

\section*{Introduction}

\bookmark[page=2,level=1]{Introduction}

Rough paths provide a rich theory of integration, beginning with Lyons' seminal work \cite{Lyo98} featuring the integration of one-forms with respect to geometric rough paths. Gubinelli then introduced a class of controlled paths, for which a theory of integration is developed in \cite{Gub04} for the $2 \leq p < 3$ case, and is extended to a more general framework of branched rough paths as the driving paths in \cite{Gub10}. These works and subsequent ones provide a solid core foundation for rough integration theory in one time variable.

As for the multivariable case, previous works such as \cite{CG14} and \cite{GH19} have succeeded in defining two-parameter rough integrals in the $2 \leq p < 3$ case. In the approach of \cite{GH19}, a class of jointly controlled paths is chosen to be such that they are twice rough integrable as controlled paths, and a Fubini type theorem is established under the assumption of admitting a smooth approximation. Given these double rough integrals, a natural question to ask would be whether we could construct a two-parameter rough integral which serves as an analogue to the integral with respect to a product measure in the classical Fubini's theorem.

Here we look to develop the work on rough Fubini type theorems done in \cite{GH19} in several ways. Firstly to extend the notion of jointly controlled paths to consider geometric driving rough paths of arbitrary $p$-variation regularity. Secondly to provide a definition of a two-parameter rough integral of a jointly controlled path, corresponding to the integral against a product measure in the classical case. We give conditions under which this integral is well-defined and establish bounds on the integral which naturally extend those obtained in the one-parameter setting. Using this new integral, which we will call the joint integral, we are able to establish a Fubini type theorem in the geometric arbitrary $p$-variation case. In the process, we are also able to relax the assumption that the integrand is a smoothly approximable path which has been made in the preceding work \cite{GH19}.

We begin Section 1 by defining jointly controlled paths in two variables with respect to one driving geometric rough path. The basic idea behind these jointly controlled paths is that the two-parameter path defines a family of controlled paths in each time variable separately and that the Gubinelli derivatives are also, in some sense, controlled paths satisfying some symmetry constraints. Just as in \cite{GH19}, we demonstrate that these jointly controlled paths can be integrated twice with respect to the controlling underlying rough path.

Once the class of jointly controlled paths are suitably defined and their double integrals shown to exist, we move on to construct the joint integral in Sections 2 and 3. Our approach to the construction of this integral is akin to the use of Young's argument in the construction of the one-parameter rough integral \cite{You36}, whereby we first bound the discrete integrals with a maximal inequality and leverage this bound to show existence and uniqueness of the limit of discrete integrals as their partitions become finer.

In Section 2 we derive a maximal inequality in the context of jointly controlled paths, heavily inspired by the arguments in the two-parameter Young-Towghi maximal inequality of \cite{Tow02} as it is presented in \cite{FV11}. Given a discrete integral over an arbitrary partition, we carefully choose a point to remove from the partition such that changes to the discrete integral are small, and then repeatedly remove points in this manner to attain a bound independent of the partition. Within this section we also establish some important intermediate bounds for future use.

Using the maximal inequality and intermediate bounds, existence of the joint integral is established in Section 3 under some mild uniformity and regularity conditions on mixed variation, a two-parameter analogue to $p$-variation used in \cite{Tow02}. Here the arguments follow those of the classical sewing lemma \cite{FdLP06}, but significant modifications are needed to adapt it to the two-parameter setting. Similar arguments and estimates to those used in showing the existence of the joint integral are then used to prove a Fubini type theorem for rough integration, equating the joint integral to the two double integrals. Following arguments from Theorem 4.17 of \cite{FH14} we also establish a stability property of the double integrals, by which we mean that a "distance" can be described such that two jointly controlled paths which are close to each other in this sense are such that their double integrals are also close.

We then generalise to the case of two different controlling rough paths in Section 5, where jointly controlled paths are suitably defined and the key results of the preceding sections are given in this context. To finish, we give an example of a jointly controlled path in the form of the signature kernel, an object which also satisfies a two-parameter rough integral equation and has seen recent use in data science applications in \cite{sigkernelsgoursat}.

\section{Jointly controlled paths}

We will assume that the reader is familiar with the basic definitions and properties of rough paths and refer to \cite{Lyo98, CLL07, FH14} as comprehensive accounts of the theory. For the reader's convenience, some key results on the integration of controlled paths are contained in the appendix, as well as notation choices made throughout.

In order to study two-parameter rough integrals, we first need to find a class of two-parameter paths which are twice integrable under one-parameter rough integration. In the recent work \cite{GH19} this is done by introducing the notion of jointly controlled paths, for which they show that the one-parameter rough integrals of these jointly controlled paths are themselves controlled paths. An alternative method of constructing double rough integrals using "rough sheets" is done in \cite{CG14}. Here we will adopt the approach of section 5 of \cite{GH19} which handles the case where $X$ is a $p$-rough path for $2 \leq p < 3$ with H\"older control. Here we extend these notions and arguments to the case where $X$ is a geometric rough path of arbitrary $p$-variation and arbitrary control. The following definition is a natural extension of Definition 5.1 in \cite{GH19}.

\begin{definition}[Jointly controlled two-parameter paths]
Let $p > 1$, $X = (1, X^1, X^2, \dotsc, X^{\pfloor})$ be an $\omega$-controlled geometric $p$-rough path on a Banach space $V$, and let $N = \pfloor - 1$. A two parameter path $Y: [0,T] \times [0,T] \rightarrow E$ on a Banach space $E$ is jointly controlled by the rough path $X$ if it satisfies the following conditions:
\begin{enumerate}[topsep=0pt]
    \item Let $Y^{(1;0,0)}_{s, \cdot} = Y_{s, \cdot}$ and $Y^{(2;0,0)}_{\cdot, u}= Y_{\cdot, u}$. For $j,k = 0, \dotsc, N$ and every $s, u \in [0,T]$, there exists
    \begin{align*}
        Y^{(1; j, k)}_{s, \cdot}
        &:
        [0, T]
        \rightarrow
        \Hom\left(
            V^{\otimes k},
            \Hom(V^{\otimes j}, E)
        \right)
        \\
        Y^{(2; k, j)}_{\cdot, u}
        &:
        [0, T]
        \rightarrow
        \Hom\left(
            V^{\otimes j},
            \Hom(V^{\otimes k}, E)
        \right)
    \end{align*}
    such that the tuples $\left( Y^{(1;j,0)}_{s, \cdot}, Y^{(1;j,1)}_{s, \cdot}, \dotsc, Y^{(1;j, N)}_{s, \cdot} \right)$ and $\left( Y^{(2;k,0)}_{\cdot, u}, Y^{(2;k,1)}_{\cdot, u}, \dotsc, Y^{(2;k, N)}_{\cdot, u} \right)$ are $X$-controlled paths (see Definition \ref{controlled path definition}).
    
    \item For all $j,k = 0, 1, \dotsc, N$, and $s, u \in [0,T]$, the derivatives $Y^{(i;j,k)}$ satisfy the symmetry condition
    \begin{equation}\label{derivative symmetry}
        Y^{(1;j,k)}_{s,u}(y)(x)
        =
        Y^{(2;k,j)}_{s,u}(x)(y)
    \end{equation}
    for $x \in V^{\otimes j}$ and $y \in V^{\otimes k}$.
\end{enumerate}
The collection $\left\{Y^{(i; j, k)} \mid i = 1, 2; \; j, k = 0 \dotsc, N\right\}$ then defines this two-parameter jointly controlled path over $[0,T]^2$. Denote by $\Dscr_X^p \left( [0,T]^2; E \right)$ this class of jointly $X$-controlled paths.
\end{definition}

From here on we will take $X$ to be a geometric $p$-rough path on a finite-dimensional Banach space $V$ unless explicitly stated otherwise. Under the symmetry condition we have that the Gubinelli derivatives $Y^{(1;0,j)}_{s, \cdot}$ of $Y_{s, \cdot}$ are such that $Y^{(1;0,j)}_{s,u} = Y^{(2;j,0)}_{s,u}$ and similarly for the derivatives of $Y_{\cdot, u}$. As such we can see that this definition arises from the idea that the Gubinelli derivatives of $Y_{s, \cdot}$ and $Y_{\cdot, u}$ are themselves controlled paths. We will use $Y$ to refer to both the collection $\{Y^{(i;j,k)}\}$ and the base path $Y^{(i;0,0)}$, depending on the context. Following from the definition of $X$-controlled paths, for any $\{Y^{(i;j,k)}\} \in \Dscr_X^p([0,T]^2; E)$ we have remainders $R^{(i; j, k)}$ defined by
\begin{align}
    \label{remainder 1}
    Y_{s, v}^{(1;j,k)}
    &=
    \sum^{N-k}_{l = 0}
        Y_{s, u}^{(1;j,k+l)}
        \left(
            X^l_{u,v}
        \right)
    +
    R^{(1;j,k)}_{s;u,v}, \\[0.5em]
    \label{remainder 2}
    Y_{t, u}^{(2;k,j)}
    &=
    \sum^{N-j}_{m = 0}
        Y_{s, u}^{(2;k, j+m)}
        \left(
            X^m_{s,t}
        \right)
    +
    R^{(2;k,j)}_{u;s,t},
\end{align}
which satisfy $\| R^{(i;j,k)}_{u} \|_{(\pfloor - k)/p} < \infty$ for any $u \in [0,T]$. Just as in the $2 \leq p < 3$ case in \cite{GH19}, the remainders satisfy some sort of symmetry relation.

\begin{lemma}\label{remainder relation lemma}
Let $\{Y^{(i;j,k)}\} \in \Dscr^p_X([0,T]^2; E)$. Let the remainders $R^{(i;j,k)}$ be defined by \eqref{remainder 1} and \eqref{remainder 2}. Define the maps $\Rbf^{(1; j, k)} : [0,T]^2 \times [0,T]^2 \rightarrow \Hom(V^{\otimes k}, \Hom(V^{\otimes j}, E))$ and $\Rbf^{(2;k, j)} : [0,T]^2 \times [0,T]^2 \rightarrow \Hom(V^{\otimes j}, \Hom(V^{\otimes k}, E))$ by
\begin{align*}
    \Rbf^{(1; j, k)}_{s,t;u,v}
    (y)(x)
    &=
    R^{(1;j,k)}_{t;u,v}(y)(x)
    -
    R^{(1;j,k)}_{s;u,v}(y)(x)
    -
    \sum^{N-j}_{m=1}
        R^{(1;j+m,k)}_{s;u,v}
        (y)
        \left(
            X^m_{s,t}
        \right)
        (x), \\[0.5em]
    \Rbf^{(2; k, j)}_{u,v;s,t}
    (x)(y)
    &=
    R^{(2;k,j)}_{v;s,t}(x)(y)
    -
    R^{(2;k,j)}_{u;s,t}(x)(y)
    -
    \sum^{N-k}_{l=1}
        R^{(2;k+l,j)}_{u;s,t}
        (x)
        \left(
            X^l_{u,v}
        \right)
        (y),
\end{align*}
where $x \in V^{\otimes j}$, $y \in V^{\otimes k}$. The remainders satisfy
\begin{equation}
    \label{remainder relation}
    \Rbf^{(1; j, k)}_{s,t;u,v}
    (y)(x)
    =
    \Rbf^{(2; k, j)}_{u,v;s,t}
    (x)(y).
\end{equation}
\end{lemma}

\begin{proof}
The equality \eqref{remainder relation} follows from simple calculations using equations \eqref{remainder 1} and \eqref{remainder 2}:
\begin{align*}
    R^{(1;j,k)}_{t;u,v}(y)(x)
    -
    R^{(1;j,k)}_{s;u,v}(y)(x)
    &=
    \left(
        Y^{(1;j,k)}_{t,v} - Y^{(1;j,k)}_{t,u}
        - Y^{(1;j,k)}_{s,v} + Y^{(1;j,k)}_{s,u}
    \right)(y)(x) \\[0.5em]
    &\quad -
    \sum^{N-k}_{l = 1}
    \left(
        Y^{(1;j,k+l)}_{t,u}
        -
        Y^{(1;j,k+l)}_{s,u}
    \right)
    \left(
        X^l_{u,v}
    \right)
    (y)(x) \\[0.5em]
    &=
    \left(
        Y^{(1;j,k)}_{t,v} - Y^{(1;j,k)}_{t,u}
        - Y^{(1;j,k)}_{s,v} + Y^{(1;j,k)}_{s,u}
    \right)(y)(x) \\[0.5em]
    &\quad -
    \sum^{N-k}_{l = 1}
        \sum^{N-j}_{m = 1}
            Y^{(2;k+l,j+m)}_{s,u}
        \left(
            X^m_{s,t}
        \right)
    (x)
    \left(
        X^l_{u,v}
    \right)
    (y) \\[0.5em]
    &\quad -
    \sum^{N-k}_{l = 1}
        R^{(2;k+l, j)}_{u;s,t}
        (x)
    \left(
        X^l_{u,v}
    \right)
    (y).
\end{align*}
A similar calculation then yields
\begin{align*}
    R^{(2;k,j)}_{v;s,t}(x)(y)
    -
    R^{(2;k,j)}_{u;s,t}(x)(y)
        &=
    \left(
        Y^{(1;j,k)}_{t,v} - Y^{(1;j,k)}_{t,u}
        - Y^{(1;j,k)}_{s,v} + Y^{(1;j,k)}_{s,u}
    \right)(y)(x)
    \\[0.5em]
    &\quad -
    \sum^{N-j}_{m = 1}
        \sum^{N-k}_{l = 1}
            Y^{(1;j+m,k+l)}_{s,u}
        \left(
            X^l_{u,v}
        \right)
    (y)
    \left(
        X^m_{s,t}
    \right)
    (x) \\[0.5em]
    &\quad -
    \sum^{N-j}_{m = 1}
        R^{(1;j+m,k)}_{s;u,v}
        (y)
    \left(
        X^m_{s,t}
    \right)
    (x) \\[0.5em]
    &=
    R^{(1;j,k)}_{t;u,v}(y)(x)
    -
    R^{(1;j,k)}_{s;u,v}(y)(x)
    +
    \sum^{N-k}_{l = 1}
        R^{(2;k+l, j)}_{u;s,t}
        (x)
    \left(
        X^l_{u,v}
    \right)
    (y) \\[0.5em]
    &\quad -
    \sum^{N-j}_{m = 1}
        R^{(1;j+m,k)}_{s;u,v}
        (y)
    \left(
        X^m_{s,t}
    \right)
    (x).
\end{align*}
Rearranging the terms then gives equation \eqref{remainder relation} as required.
\end{proof}

Again, similarly to the $2 \leq p < 3$ case in \cite{GH19}, this relation allows us to show that the remainders themselves are $X$-controlled paths in some sense.

\begin{corollary}\label{remainders are controlled paths}
For $i = 1, 2$, $j,k = 0, \dotsc, N$, define the map $R^{(i;j,k)}: [0, T]^2 \times [0, T] \rightarrow \Hom(V^{\otimes k}, \Hom(V^{\otimes j}, E))$ by
\[
    R^{(i;j,k)}_{s,t;u}(x)(y)
    =
    R^{(i';k,j)}_{u;s,t}(y)(x)
\]
for $x \in V^{\otimes j}$ and $y \in V^{\otimes k}$, where $i' \neq i$. For fixed $0 \leq s < t \leq T$, the tuple $\left( R^{(i;j,k)}_{s,t}, \dotsc R^{(i;j,N)}_{s,t} \right)$ is an $X$-controlled path with remainders $\Rbf^{(i;j,k)}_{s,t} : [0, T]^2 \rightarrow \Hom(V^{\otimes k}, \Hom(V^{\otimes j}, E))$ as defined in Lemma \ref{remainder relation lemma}.
\end{corollary}

\begin{proof}
By Lemma \ref{remainder relation lemma} we know that
\[
    R^{(i;j,k)}_{s,t;v}
    =
    \sum^{N-k}_{l=0}
        R^{(i;j, k+l)}_{s,t;u}
        \left(
            X^l_{u,v}
        \right)
    +
    \Rbf^{(i;j,k)}_{s,t;u,v}
\]
and from the definition of $\Rbf^{(i;j,k)}$ it is easy to show that $\| \Rbf^{(i;j,k)}_{s,t;u,v} \|_{p/(\pfloor - k)}$ is finite for any fixed choice of $s,t$ with $0 \leq s \leq t \leq T$.
\end{proof}

One quantity which we will need to consider is the mixed variation introduced in \cite{Tow02}, which also is used in works such as \cite{FGGR16}. It is the two-parameter analogue to $p$-variation.

\begin{definition}
Given $p, q  \geq 1$ and a process $A:[0, T]^2 \times [0,T]^2 \rightarrow E$ on a Banach space $E$ define the mixed $(p,q)$-variation over $[s, t]\times[u,v]$ by
\[
    V^{p,q}_{[s,t]\times[u,v]}(A)
    :=
    \sup_{\Dcal \times \Dcal'}
    \left(
        \sum^{n_0-1}_{n=0}
        \left(
            \sum^{m_0-1}_{m=0}
            \left|
                A
                \binom{s_m, s_{m+1}}{u_n, u_{n+1}}
            \right|^{p}
        \right)^{q/p}
    \right)^{1/q}
\]
where the supremum is over partitions $\Dcal = \{ s_0 < \dotsc < s_{m_0}\} \subset [s,t]$ and $\Dcal' = \{ u_0 < \dotsc < u_{n_0}\} \subset [u, v]$. If $V^{p,q}_{[s,t]\times[u,v]}(A)$ is finite then we say that $A$ has finite $(p,q)$-variation over $[s,t] \times [u,v]$.

We will also say that given a one-dimensional control $\omega$ that $A$ has finite $\omega$-controlled $(p,q)$-variation over $[0,T]^2$ if $\| A \|_{p,q,\omega}$ is finite, where we define the norm $\| \cdot \|_{p,q,\omega}$ by
\[
    \| A \|_{p,q,\omega}
    :=
    \sup_{\substack{
            0 \leq s < t \leq T
            \\
            0 \leq u < v \leq T
        }
    }
    \left|
        \frac{
            A\binom{s,t}{u,v}
        }{
            \omega(s,t)^{1/p}
            \omega(u,v)^{1/q}
        }
    \right|
    .
\]
\end{definition}
By subadditivity of controls, finite $\omega$-controlled $(p,q)$-variation implies finite mixed $(p,q)$-variation. With this in mind we place a mixed variation condition on the second order remainders $\Rbf^{(i;j,k)}$.

\begin{condition} \label{mixed variation condition}
Let $\{Y^{(i;j,k)}\} \in \Dscr^p_X([0,T]^2; E)$. For $l = 0, \dotsc, N$, let $p_l:= p/(l+1)$ and suppose that there exists $q_l \geq 1$ such that $\frac{1}{p_l} + \frac{1}{q_l} =: \theta_l > 1$. Assume for $j,k = 0 \dotsc, N$ that $R^{(i;j,k)}_u$ has finite $\omega$-controlled $q_k$-variation and  $\Rbf^{(1;j,k)}$ have finite $(q_{j},q_{k})$-variation, where we write
\[
    \Rbf^{(i;j,k)}
    \binom{s,t}{u,v}
    =
    \Rbf^{(i;j,k)}_{s,t;u,v}
    .
\]
Denote by $\theta_* = \min_{k}\theta_k$ and $\theta^* = \max_{k}\theta_k$.
\end{condition}

Taking $q_l = \frac{p}{\pfloor - l}$, by definition we know that the remainders $R^{(i;j,k)}_u$ have finite $\omega$-controlled $q_k$-variation and thus $\Rbf^{(i;j,k)}_{s,t}$ has finite $q_k$-variation. By the remainder relation \eqref{remainder relation} it then follows that $\Rbf^{(i;j,k)}_{\cdot, \cdot; u,v}$ has finite $q_j$-variation.

Since $\left(Y_{r, \cdot} = Y^{(1;0,0)}_{r, \cdot}, \dotsc, Y^{(1;0,N)}_{r, \cdot}\right)$ and $\left(Y_{\cdot, r'} = Y^{(2;0,0)}_{\cdot, r'}, \dotsc, Y^{(2;0,N)}_{\cdot, r'} \right)$ are controlled paths, we know that the integrals
\begin{alignat}{4}
    \label{integral 1}
    Z^{(1)}_r
    &=
    \int_u^v
        Y_{r,r'}
    \,dX_{r'}
    &&=
    \lim_{|\Dcal'| \rightarrow 0}
        &&\sum_{\{u_{n-1}, u_{n}\} \subset \Dcal'} &&
            \sum_{k = 0}^{N}
                Y_{r, u_{n-1}}^{(1;0,k)}
                \left(
                    X^{k+1}_{u_{n-1}, u_{n}}
                \right)
    \\
    \label{integral 2}
    Z^{(2)}_{r'}
    &=
    \int_s^t
        Y_{r,r'}
    \,dX_{r}
    &&=
    \lim_{|\Dcal| \rightarrow 0}
        &&\sum_{\{s_{m-1}, s_{m}\} \subset \Dcal} &&
            \sum_{j = 0}^{N}
                Y_{s_{m-1}, r'}^{(2;0,j)}
                \left(
                    X^{j+1}_{s_{m+1}, s_{m}}
                \right)
\end{alignat}
exist and are well defined. In the next lemma we show that if $E = \Hom(V, W)$ for some Banach space $W$, then $Z_r^{(1)}$ and $Z_r^{(2)}$ are themselves $X$-controlled paths.

\begin{lemma}\label{integrals are controlled paths}
Let $\{Y^{i;j,k}\} \in \Dscr^p_X([0,T]^2; \Hom(V,W))$ for some Banach space $W$. Suppose that Condition \ref{mixed variation condition} is satisfied and moreover the mixed variations of $\Rbf^{(i;j,k)}$ are $\omega$-controlled. Consider intervals $[s,t]$ and $[u,v]$ in $[0, T]$, and let $Z^{(1)}, Z^{(2)}$ be defined as in \eqref{integral 1} and \eqref{integral 2}. Define also $Z^{(1;j)}:[s,t] \rightarrow \Hom(V^{\otimes j}, W)$ and $Z^{(2;k)}:[u,v] \rightarrow \Hom(V^{\otimes k}, W$) by
\begin{align*}
    Z_r^{(1;j)}
    =
    \int_u^v
        Y_{r,r'}^{(1;j,0)}
    \,dX_{r'},
    \qquad
    Z_{r'}^{(2;k)}
    =
    \int_s^t
        Y_{r,r'}^{(2;k,0)}
    \,dX_{r}
\end{align*}
for $j,k = 1, \dotsc, N$ and $Z^{(i;0)} = Z^{(i)}$. The tuples $(Z^{(1;0)}, \dotsc, Z^{(1;N)})$ and $(Z^{(2;0)}, \dotsc, Z^{(2;N)})$ are $X$-controlled paths.
\end{lemma}

\begin{remark}
Before we proceed we make some brief remarks on some upcoming notation. The local approximations of $Z_r^{(1;j)}$ are of the form
\[
    \sum_{k = 0}^{N}
        Y_{r, u}^{(1;j,k)}
        \left(
            X^{k+1}_{u,v}
        \right)
\]
where we interpret $Y_{r, r'}^{(i;j,k)}\big(X^{k+1}_{a,b}\big) \in \Hom(V^{\otimes j}, W)$ as the map sending $x \in V^{\otimes j}$ to the canonical injection of $Y_{r, r'}^{(i;j,k)}(\cdot)(x) \in \Hom(V^{\otimes k}, \Hom(V, W))$ into $\Hom(V^{\otimes k+1}, W)$ evaluated at $X^{k+1}_{a,b} \in V^{\otimes k+1}$. In this manner we may also write the local approximations of $Z^{(2;k)}_{r'}$ in the form
\[
    \sum_{j = 0}^{N}
        Y_{s, r'}^{(2;k,j)}
        \left(
            X^{j+1}_{s,t}
        \right)
    .
\]
Given a map $\Xi: [0,T]^2 \rightarrow E$ for some vector space $E$, we will use $\delta \Xi : [0,T]^3 \rightarrow E$ to be the map $\delta\Xi_{s,s',t} = \Xi_{s,t} - \Xi_{s,s'} - \Xi_{s',t}$.
\end{remark}

\begin{proof}
We first note that we may write
\begin{align*}
    Y_{t, r'}^{(1;j,k)}(y)(x)
    =
    Y_{t, r'}^{(2;k,j)}(x)(y)
    &=
    \sum^{N-j}_{m = 0}
        Y_{s, r'}^{(2;k, j+m)}
        \left(
            X^m_{s,t}
        \right)
        (x)
        (y)
    +
    R^{(2;k,j)}_{r';s,t}
    (x)(y)
    \\
    &=
    \sum^{N-j}_{m = 0}
        Y_{s, r'}^{(1;j+m,k)}
        (y)
        \left(
            X^m_{s,t}
        \right)
        (x)
    +
    R^{(1;j,k)}_{s,t;r'}(y)(x)
\end{align*}
and so for each $j = 0, \dotsc, N$ we decompose the controlled path $\left(Y^{(1;j,0)}_{t, \cdot}, \dotsc, Y^{(1;j, N)}_{t, \cdot} \right)$ into the sum of controlled paths $\left(R^{(1;j,0)}_{s,t}, \cdots, R^{(1;j,N)}_{s,t}\right)$ and $\left(Y^{(1;j+m,0)}_{s, \cdot}(\cdot)(X^m_{s,t}), \dotsc, Y^{(1;j+m, N)}_{s, \cdot}(\cdot)(X^m_{s,t}) \right)$ for $m = 0, \dotsc, N-j$. Under this decomposition we write
\begin{align*}
    Z_t^{(1;j)}
    =
    \int_u^v
        Y_{t,r'}^{(1; j)}
    \,dX_{r'}
    &=
    \int_u^v
        \left(
            \sum^{N-j}_{m = 0}
                Y_{s, r'}^{(1;j+m,0)}
                (\cdot)
                \left(
                    X^m_{s,t}
                \right)
            +
            R^{(1;j,0)}_{s,t;r'}
        \right)
    \,dX_{r'}
    \\[0.5em]
    &=
    \sum^{N-j}_{m = 0}
        Z_s^{(j+m)}
        \left(
            X^m_{s,t}
        \right)
    +
    \int_u^v
        R^{(1;j,0)}_{s,t;r'}
    \,dX_{r'}
\end{align*}
We now look at the regularity of the remainder term above, denote said term by $R^{(Z;1;j)}_{s,t}$. Fixing $r' \in [0,T]$ and writing
\[
    \Xi_{s,t;u,v}^{(j)}
    =
    \sum^N_{k=0}
        R^{(1;j,k)}_{s,t;u}
        \left(
            X^{k+1}_{u,v}
        \right)
    ,
\]
from Lemma \ref{controlled path identity lemma} and Corollary \ref{remainders are controlled paths} we deduce
\[
    \left\|
        \delta \Xi^{(j)}_{s,t}
    \right\|_{1/\theta_*}
    \leq
    C(u,v)
    \sum^N_{k=0}
        \left\|
            \Rbf^{(1;j,k)}_{s,t}
        \right\|_{q_k}
        \left\|
            X^{k+1}
        \right\|_{p_k}
    .
\]
Applying Lemma \ref{sewing lemma} we then have
\begin{align*}
    \left|
        R^{(Z;1;j)}_{s,t}
    \right|
    &\leq
    \left|
        \Xi_{s,t;u,v}^{(j)}
    \right|
    +
    \zeta(\theta_*)
    \;
    \omega(u,v)^{\theta_*}
    \left\|
        \delta \Xi^{(j)}_{s,t}
    \right\|_{1/\theta_*}
    \\[0.5em]
    &\leq
    \sum^N_{k=0}
        \left\|
            R^{(2;k,j)}_u
        \right\|_{q_j}
        \left|
            X^{k+1}_{u,v}
        \right|
        \omega(s,t)^{1/q_j}
    \\
    & \quad
    +
    C(u,v)
    \;
    \zeta(\theta_*)
    \;
    \omega(u,v)^{\theta_*}
    \sum^N_{k=0}
        \left\|
            \Rbf^{(1;j,k)}_{s,t}
        \right\|_{q_j,q_k}
        \left\|
            X^{k+1}
        \right\|_{p_k}
    ,
\end{align*}
which tells us that $\|R^{(Z;1;j)}\|_{q_j}$ is finite for $j=0, \dotsc, N$, and thus the tuple $(Z^{(1;0)}, \dotsc, Z^{(1;N)})$ is an $X$-controlled path. The proof for showing that $(Z^{(2;0)}, \dotsc, Z^{(2;N)})$ is a controlled path follows almost identically.
\end{proof}

Setting the space $W= \Hom(V, U)$ for some Banach space $U$, we can then define the double rough integrals
\begin{alignat}{2}
    \Ical_1
    \binom{s,t}{u,v}
    &:=
    \int^t_s
        \left(
            \int^v_u
                Y_{r,r'}
            \,dX_{r'}
        \right)
    \,dX_r
    &&=
    \int^t_s
        Z^{(1)}_r
    \,dX_r,
    \\[0.5em]
    \Ical_2
    \binom{u,v}{s,t}
    &:=
    \int^v_u
        \left(
            \int^t_s
                Y_{r,r'}
            \,dX_{r}
        \right)
    \,dX_{r'}
    &&=
    \int^t_s
        Z^{(2)}_{r'}
    \,dX_{r'}.
\end{alignat}

With these double rough integrals constructed, we now look to show a Fubini's theorem for these integrals. To this end, we consider a third integral which is analogous to the integral with respect to a product measure in the classical Fubini's theorem. If we naively approximate these integrals twice, we observe that the two integrals share the same local approximation due to the symmetry condition on our jointly controlled paths,
\begin{equation}
    \sum_{j, k = 0}^{N}
        Y_{s,u}^{(1;j,k)}
        \Big(
            X^{k+1}_{u,v}
        \Big)
        \Big(
            X^{j+1}_{s,t}
        \Big)
    =
    \sum_{j, k = 0}^{N}
        Y_{s,u}^{(2;k,j)}
        \Big(
            X^{j+1}_{s,t}
        \Big)
        \Big(
            X^{k+1}_{u,v}
        \Big)
\end{equation}
where we understand $Y_{r,r'}^{(i;j,k)}\big(X^{k+1}_{a,b}\big)\big(X^{j+1}_{c,d}\big)$ by taking $Y_{r,r'}^{(i;j,k)}\big(X^{k+1}_{a,b}\big)$ as before and identifying it with its canonical injection into $\Hom(V^{\otimes j+1}, U)$. With this local approximation as a starting point we move towards construction of a third integral through a Young type argument, with inspiration taken from the two-parameter Young-Towghi inequality of \cite{Tow02} as presented in \cite{FV11}.

\section{A Young type maximal inequality}

We now begin our first steps towards a Fubini type theorem for rough integration of jointly controlled paths. Our approach is akin to the case of Young integration, whereby we first establish a "maximal inequality" over the discrete two-parameter integrals and then leverage this inequality to prove existence and uniqueness of the integral as the limit of the discrete two-parameter integrals over partitions with decreasing mesh size. A Fubini type theorem then follows from similar estimates to those used in proving the existence of the integral. We are also inspired by the sewing lemma of \cite{FdLP06}, as can be seen in Lemma \ref{sewing lemma} in the appendix, where the proof of the sewing lemma is formatted to mirror the same type of argument used here.

In comparison with \cite{GH19} which also establishes a Fubini type theorem for rough integration, we are able to generalise to geometric driving rough paths of arbitrary $p$-variation and also remove the restriction of being smoothly approximable. The result to be given here also draws more parallels with the classical Fubini's theorem by providing a third integral constructed directly from the local approximations, serving as an analogue to the integral with respect to a product measure. We will refer to the iterated rough integrals as double integrals and refer to the third type of integral as the joint rough integral.

Borrowing from \cite{FV11}, we will use grid-like partitions of $[s,t] \times [u,v]$ to refer to partitions of the form $\Dcal \times \Dcal'$ where $\Dcal$ and $\Dcal'$ are partitions of $[s,t]$ and $[u,v]$ respectively. For a function $\Omega : [0, T]^2 \times [0, T]^2 \rightarrow E$ on a Banach space $E$ and partitions $\Dcal = \{ s_0 < \dotsc < s_{m_0}\} \subset [s,t]$ and $\Dcal' = \{ u_0 < \dotsc < u_{n_0} \} \subset [u,v]$ we write
\[
    \sum_{\Dcal \times \Dcal'}
        \Omega
    :=
    \sum_{m=1}^{m_0}
    \sum_{n=1}^{n_0}
        \Omega
        \binom{s_{m}, s_{m+1}}{u_{n}, u_{n+1}}
    ,
\]
which will be used to write the discrete two-parameter integral over $\Dcal \times \Dcal'$ in the context of two parameter rough integration.

We follow the same type of argument as in the two dimensional Young-Towghi maximal inequality of \cite{Tow02} as presented in the appendix of \cite{FV11}. The main idea is as follows: given a discrete integral over a partition, we select a point to remove from the partition and then observe the change in the discrete integral. The particular point we remove is carefully chosen to keep this change small. We then repeat until we are left with the trivial partition. To keep notation succinct, we introduce the following quantities.

\begin{definition}
\label{local approximation definitions}
Let $\{Y^{(i;j,k)}\}$ be a jointly $X$-controlled two-parameter path on $\Hom(V, V^*)$. Denote by $\Omega^Y: [0,T]^2 \times [0,T]^2 \rightarrow \Real$ the local approximation of the joint integral
\[
    \Omega^Y
    \binom{s,t}{u,v}
    =
    \sum_{j,k=0}^N
        Y_{s,u}^{(1;j,k)}
        \Big(
            X^{k+1}_{u,v}
        \Big)
        \Big(
            X^{j+1}_{s,t}
        \Big)
    =
    \sum_{j,k=0}^N
        Y_{s,u}^{(2;j,k)}
        \Big(
            X^{k+1}_{s,t}
        \Big)
        \Big(
            X^{j+1}_{u,v}
        \Big).
\]
Now define $\Gamma^Y: [0,T]^3 \times [0,T]^2 \rightarrow \Real$ by
\begin{equation*}
    \Gamma^Y
    \binom{s,s',t}{u,v}
    :=
    -\delta\left(
        \Omega^Y
        \binom{\cdot,\cdot}{u,v}
    \right)
    \left(
        s,s',t
    \right)
    =
    \Omega^Y
    \binom{s,s'}{u,v}
    +
    \Omega^Y
    \binom{s',t}{u,v}
    -
    \Omega^Y
    \binom{s,t}{u,v}
\end{equation*}
and $\Theta^Y: [0,T]^3 \times [0,T]^3 \rightarrow \Real$ by
\begin{equation*}
    \Theta^Y
    \binom{s,s',t}{u,u',v}
    :=
    -\delta\left(
        \Gamma^Y
        \binom{s,s',t}{\cdot,\cdot}
    \right)
    \left(
        u,u',v
    \right)
    =
    \Gamma^Y
    \binom{s,s',t}{u,u'}
    +
    \Gamma^Y
    \binom{s,s',t}{u',v}
    -
    \Gamma^Y
    \binom{s,s',t}{u,v},
\end{equation*}
where we recall that for $\Xi : [0, T]^2 \rightarrow E$ that $\delta \Xi_{s,s',t} = \Xi_{s,t} - \Xi_{s,s'} - \Xi_{s',t}$ for $0 \leq s < s' < t \leq T$.
\end{definition}

\begin{remark}
We may also similarly define $\Gamma^Y: [0,T]^2 \times [0,T]^3 \rightarrow \Real$ by instead fixing the first pair of variables,
\[
    \Gamma^Y
    \binom{s,t}{u,u',v}
    :=
    -\delta\left(
        \Omega^Y
        \binom{s,t}{\cdot,\cdot}
    \right)
    \left(
        u,u',v
    \right)
\]
and it is easily verified that we may alternatively define $\Theta^Y$ by
\[
    \Theta^Y
    \binom{s,s',t}{u,u',v}
    :=
    -\delta\left(
        \Gamma^Y
        \binom{\cdot,\cdot}{u,u',v}
    \right)
    \left(
        s,s',t
    \right)
\]
due to the symmetry of $\Omega^Y$.
\end{remark}

From Lemma \ref{controlled path identity lemma} we show that $\Gamma^Y$ and $\Theta^Y$ can be expressed in terms of first and second order remainders of $Y$.

\begin{lemma}
Let $\{Y^{(i;j,k)}\}$ be a jointly $X$-controlled two-parameter path on $\Hom(V,V^*)$. Then the following identities hold:
\begin{align}
    \label{Gamma identity}
    \Gamma^Y
    \binom{s,s',t}{u,v}
    &=
    \sum_{j,k=0}^N
        R^{(2;k,j)}_{u;s,s'}
        \Big(
            X^{j+1}_{s',t}
        \Big)
        \Big(
            X^{k+1}_{u,v}
        \Big),
    \\[0.5em]
    \label{Theta identity}
    \Theta^Y
    \binom{s,s',t}{u,u',v}
    &=
    \sum_{j,k=0}^N
        \Rbf^{(1;j,k)}_{s,s';u,u'}
        \Big(
            X^{k+1}_{u',v}
        \Big)
        \Big(
            X^{j+1}_{s',t}
        \Big)
    .
\end{align}
\end{lemma}

\begin{proof}
By definition we have
\begin{align*}
    \Gamma^Y
    \binom{s,s',t}{u,v}
    &=
    \sum^N_{j,k=0}
        Y^{(1;j,k)}_{s,u}
        \Big(
            X^{k+1}_{u,v}
        \Big)
        \Big(
            X^{j+1}_{s,s'}
        \Big)
        +
        Y^{(1;j,k)}_{s',u}
        \Big(
            X^{k+1}_{u,v}
        \Big)
        \Big(
            X^{j+1}_{s',t}
        \Big)
        -
        Y^{(1;j,k)}_{s,u}
        \Big(
            X^{k+1}_{u,v}
        \Big)
        \Big(
            X^{j+1}_{s,t}
        \Big)
    \\[0.5em]
    &=
    \sum^N_{k=0}
    \left(
        \sum^N_{j=0}
            Y^{(2;k,j)}_{s,u}
            \Big(
                X^{j+1}_{s,s'}
            \Big)
            +
            Y^{(2;k,j)}_{s',u}
            \Big(
                X^{j+1}_{s',t}
            \Big)
            -
            Y^{(2;k,j)}_{s,u}
            \Big(
                X^{j+1}_{s,t}
            \Big)
    \right)
    \Big(
        X^{k+1}_{u,v}
    \Big).
\end{align*}
Applying Lemma \ref{controlled path identity lemma} we then arrive at \eqref{Gamma identity}. Similarly when we use this identity in the definition of $\Theta$ we get
\begin{align*}
    \Theta^Y
    \binom{s,s',t}{u,v}
    &=
    \sum^N_{j,k=0}
        R^{(2;k,j)}_{u;s,s'}
        \Big(
            X^{j+1}_{s',t}
        \Big)
        \Big(
            X^{k+1}_{u,u'}
        \Big)
        +
        R^{(2;k,j)}_{u';s,s'}
        \Big(
            X^{j+1}_{s',t}
        \Big)
        \Big(
            X^{k+1}_{u',v}
        \Big)
        -
        R^{(2;k,j)}_{u;s,s'}
        \Big(
            X^{j+1}_{s',t}
        \Big)
        \Big(
            X^{k+1}_{u,v}
        \Big)
    \\[0.5em]
    &=
    \sum^N_{j=0}
    \left(
        \sum^N_{k=0}
            R^{(1;j,k)}_{s,s';u}
            \Big(
                X^{k+1}_{u,u'}
            \Big)
            +
            R^{(1;j,k)}_{s,s';u'}
            \Big(
                X^{k+1}_{u',v}
            \Big)
            -
            R^{(1;j,k)}_{s,s';u}
            \Big(
                X^{k+1}_{u,v}
            \Big)
    \right)
    \Big(
        X^{j+1}_{s',t}
    \Big).
\end{align*}
We recall now that Lemma \ref{remainders are controlled paths} states that for fixed $s,s' \in [0,T]$ the tuples $\left( R^{(1;j,0)}_{s,s'}, \dotsc, R^{(1;j,N)}_{s,s'}\right)$ are controlled paths with remainders $\left\{ \Rbf^{(1;j,k)}_{s,s';\cdot,\cdot} \right\}_{k=0}^N$, and so by applying Lemma \ref{controlled path identity lemma} again we arrive at \eqref{Theta identity}.
\end{proof}

Now we work towards building a series of intermediate bounds, which will prove useful in proving both the maximal inequality and the rough Fubini type theorem. The cornerstone for these bounds is the following key estimate, which we make under assumptions of finite mixed variation.

\begin{lemma}\label{Theta inequality lemma}
Let $\{Y^{(i;j,k)}\}$ be a jointly $X$-controlled two-parameter path on $\Hom(V, V^*)$ and suppose that Condition \ref{mixed variation condition} holds and let $\alpha \in (1/\theta_*, 1)$. Consider partitions $\Dcal = \{s_{0}< \dotsc <s_{m_0}\} \subset [s,t]$ and $\Dcal' = \{u_{0} < \dotsc < u_{n_0} \} \subset [u,v]$, and for $n=1, \dotsc, n_0 - 1$ define $\Pcal^{(1; \Dcal, \Dcal')}_{\alpha, n}$ by
\[
    \Pcal^{(1; \Dcal, \Dcal')}_{\alpha, n}
    :=
    \sum^{m_0 -1}_{m = 1}
        \left|
            \Theta^Y
            \binom{s_{m-1}, s_m, s_{m+1}}{u_{n-1}, u_n, u_{n+1}}
        \right|^{\alpha}
\]
Then there exists $n^*$ such that 
\begin{equation}\label{variation inequality}
    \Pcal^{(1;\Dcal, \Dcal')}_{\alpha,n^*}
    \leq
    \frac{C(p,\alpha,\theta_*)}{(n_0 -1)^{\alpha \theta_*}}
    \Vbf^{(1;p)}_{[s,t],[u,v]}(X, Y),
\end{equation}
where $C(p,\alpha,\theta_*)$ is such that $C(p,\alpha,\theta_*) \leq \pfloor^{2+\alpha \theta^*}$ and
\[
    \Vbf^{(1;p)}_{[s,t],[u,v]}(X, Y)
    :=
    \max_{j,k = 0, \dotsc, N}
        \left\{
            V^{p/(j+1)}_{[s,t]}\left(
                X^{j+1}
            \right)
            V^{p/(k+1)}_{[u,v]}\left(
                X^{k+1}
            \right)
            V^{q_j, q_k}_{[s,t]\times[u,v]}\left(
                \Rbf^{(1; j, k)}
            \right)
        \right\}.
\]
\end{lemma}

\begin{proof}
Recall that $N = \pfloor - 1$, $p_j = p/(j+1)$ and $\frac{1}{p_j} + \frac{1}{q_j} = \theta_j \geq \theta_* > 1$. Since $\alpha < 1$ we have $\left(\sum_i |a_i|\right)^{\alpha} \leq \sum_i |a_i|^{\alpha}$, which gives us
\begin{align}
    \nonumber
    \Pcal^{(1; \Dcal,\Dcal')}_{\alpha, n}
    &\leq
    \sum^N_{j,k=0}
    \sum^{m_0 -1}_{m=1}
    \left|
        \Rbf^{(1; j, k)}_{s_{m-1}, s_m; u_{n-1}, u_n}
    \right|^{\alpha}
    \left|
        X^{j+1}_{s_m, s_{m+1}}
    \right|^{\alpha}
    \left|
        X^{k+1}_{u_n, u_{n+1}}
    \right|^{\alpha}
    \\[0.5em]
    \nonumber
    &\leq
    \sum^{N}_{j,k=0}
        \left(
            \sum^{m_0 -1}_{m=1}
            \left|
                \Rbf^{(1; j, k)}_{s_{m-1}, s_m; u_{n-1}, u_n}
            \right|^{\alpha \theta_j q_j}
        \right)^{1/\theta_j q_j}
        \left(
            \sum^{m_0 -1}_{m=1}
            \left|
                X^{j+1}_{s_m, s_{m+1}}
            \right|^{\alpha \theta_j p_j}
        \right)^{1/\theta_j p_j}
        \left|
            X^{k+1}_{u_n, u_{n+1}}
        \right|^{\alpha}
    \\[0.5em]
    \nonumber
    &\leq
    \sum^{N}_{j,k=0}
        \left(
            \sum^{m_0 -1}_{m=1}
            \left|
                \Rbf^{(1; j, k)}_{s_{m-1}, s_m; u_{n-1}, u_n}
            \right|^{q_j}
        \right)^{\alpha/q_j}
        \left(
            \sum^{m_0 -1}_{m=1}
            \left|
                X^{j+1}_{s_m, s_{m+1}}
            \right|^{p_j}
        \right)^{\alpha/p_j}
        \left|
            X^{k+1}_{u_n, u_{n+1}}
        \right|^{\alpha}
    \\[0.5em]
    \label{Pcal inequality}
    &\leq
    \sum^{N}_{j,k=0}
        \left(
            \sum^{m_0 -1}_{m=1}
            \left|
                \Rbf^{(1; j, k)}_{s_{m-1}, s_m; u_{n-1}, u_n}
            \right|^{q_j}
        \right)^{\alpha/q_j}
        V^{p_j}_{[s,t]}
        \left(
            X^{j+1}
        \right)^{\alpha}
        \left|
            X^{k+1}_{u_n, u_{n+1}}
        \right|^{\alpha},
\end{align}
where we use H\"older's inequality in the second line, and in the third line we use that the $\ell^{\alpha \theta_j q_j}$-norm is dominated by the $\ell^{q_j}$-norm on $\Real^{m_0 - 1}$ since $\alpha \theta_j > 1$. To shorten future notation, define
\[
    A_n^{j,k}
    :=
    \left(
        \sum^{m_0 -1}_{m=1}
        \left|
            \Rbf^{(1; j, k)}_{s_{m-1}, s_m; u_{n-1}, u_n}
        \right|^{q_j}
    \right)^{1/q_j}.
\]
For $j,k = 0 \dotsc, N$ define the sets
\[
    B_{j,k}
    :=
    \left\{
        n
        \; : \;
        \left|
            V^{p_j}_{[s,t]}
            \left(
                X^{j+1}
            \right)
            A^{j,k}_n
            X^{k+1}_{u_n, u_{n+1}}
        \right|
        =
        \max_{j',k'= 0 \dotsc, N} 
        \left|
            V^{p_{j'}}_{[s,t]}
            \left(
                X^{j'+1}
            \right)
            A^{j',k'}_n
            X^{k'+1}_{u_n, u_{n+1}}
        \right|
    \right\},
\]
which are (possibly empty) subsets of $\{1, \dotsc, n_0 - 1\}$ containing $n$ such that the pair $(j,k)$ maximises the summands in \eqref{Pcal inequality}. Since each element of $\{1, \dotsc, n_0 - 1\}$ belongs to $B_{j,k}$ for at least one pair $(j,k)$, we have that
$
    \sum_{j,k = 1}^{\pfloor}
        |B_{j,k}|
    \geq
    n_0 - 1
$
and so there must exist at least one pair $(j_0, k_0)$ such that $|B_{j_0,k_0}| \geq (n_0 - 1)/ \pfloor^2$. Let $n^* \in B_{j_0, k_0}$ be such that
\[
    \left|
        V^{p_{j_0}}_{[s,t]}
        \left(
            X^{j_0 +1}
        \right)
        A^{j_0,k_0}_{n^*}
        X^{k_0 +1}_{u_{n^*}, u_{n^*+1}}
    \right|
    =
    \min_{n \in B_{j_0, k_0}}
    \left|
        V^{p_{j_0}}_{[s,t]}
        \left(
            X^{j_0 +1}
        \right)
        A^{j_0,k_0}_n
        X^{k_0 +1}_{u_n, u_{n+1}}
    \right|.
\]
Since $n^*$ minimises the above expression over $B_{j_0, k_0}$, we obtain the following:
\begin{align*}
    \Pcal^{(1;\Dcal, \Dcal')}_{\alpha, n^*}
    &\leq
    \pfloor^2
    \left|
        V^{p_{j_0}}_{[s,t]}
        \left(
            X^{j_0 +1}
        \right)
        A^{j_0,k_0}_{n^*}
        X^{k_0 +1}_{u_{n^*}, u_{n^*+1}}
    \right|^{\alpha}
    \\[0.5em]
    &\leq
    \pfloor^2
    \left(
        V^{p_{j_0}}_{[s,t]}
        \left(
            X^{j_0 +1}
        \right)
    \right)^{\alpha}
    \left(
        \prod_{n \in B_{j_0, k_0}}
        \left|
            A^{j_0,k_0}_n
        \right|^{\alpha}
        \left|
            X^{k_0 +1}_{u_n, u_{n+1}}
        \right|^{\alpha}
    \right)^{1/|B_{j_0, k_0}|}.
\end{align*}
We now separate the product over $B_{j_0, k_0}$ above into two parts and work on them separately. By the arithmetic-geometric mean inequality we have
\begin{align*}
    \left(
        \prod_{n \in B_{j_0, k_0}}
        \left|
            A^{j_0,k_0}_n
        \right|^{q_{k_0}}
    \right)^{\alpha/q_{k_0}|B_{j_0, k_0}|}
    &\leq
    \left(
        \frac{1}{|B_{j_0, k_0}|}
        \sum_{n \in B_{j_0, k_0}}
        \left|
            A^{j_0,k_0}_n
        \right|^{q_{k_0}}
    \right)^{\alpha/q_{k_0}}
    \\[0.5em]
    & \leq
    \left(
        \frac{\pfloor^2}{n_0-1}
        \sum_{n =1}^{n_0 - 1}
        \left(
            \sum^{m_0 -1}_{m=1}
            \left|
                \Rbf^{(1; j_0, k_0)}_{s_{m-1}, s_m; u_{n-1}, u_n}
            \right|^{q_{j_0}}
        \right)^{{q_{k_0}}/{q_{j_0}}}
    \right)^{\alpha/q_{k_0}}
    \\[0.5em]
    &\leq
    \left(
        \frac{\pfloor^2}{n_0-1}
    \right)^{\alpha/q_{k_0}}
    \left(
        V^{q_{j_0}, q_{k_0}}_{[s,t]\times[u,v]}
        \left(
            \Rbf^{(1; j_0, k_0)}
        \right)
    \right)^{\alpha}.
\end{align*}
Similarly, we again apply the AM-GM inequality to get
\begin{align*}
    \left(
        \prod_{n \in B_{j_0, k_0}}
        \left|
            X^{k_0 +1}_{u_n, u_{n+1}}
        \right|^{p_{k_0}}
    \right)^{\alpha/p_{k_0}|B_{j_0, k_0}|}
    & \leq
    \left(
        \frac{1}{|B_{j_0, k_0}|}
        \sum_{n \in B_{j_0, k_0}}
        \left|
            X^{k_0 +1}_{u_n, u_{n+1}}
        \right|^{p_{k_0}}
    \right)^{\alpha/p_{k_0}}
    \\[0.5em]
    & \leq
    \left(
        \frac{\pfloor^2}{n_0-1}
    \right)^{\alpha/p_{k_0}}
    \left(
        V^{p_{k_0}}_{[u,v]}
        \left(
            X^{k_0 +1}
        \right)
    \right)^{\alpha}
\end{align*}
Combining these inequalities, we get
\begin{equation*}
    \Pcal^{(1;\Dcal, \Dcal')}_{\alpha, n^*}
    \leq
    \frac{
        \pfloor^{2 + \alpha \theta_{k_0}}
    }{
        \left(
            n_0 - 1
        \right)^{\alpha \theta_{k_0}}
    }
    \left|
        V^{p_{j_0}}_{[s,t]}
        \left(
            X^{j_0 +1}
        \right)
    \right|^{\alpha}
    \left|
        V^{p_{k_0}}_{[u,v]}
        \left(
            X^{k_0 +1}
        \right)
    \right|^{\alpha}
    \left(
        V^{q_{j_0}, q_{k_0}}_{[s,t]\times[u,v]}
        \left(
            \Rbf^{(1; j_0, k_0)}
        \right)
    \right)^{\alpha}
\end{equation*}
from which \eqref{variation inequality} follows.
\end{proof}

\begin{remark}
We will also want to consider the quantity
\[
    \Pcal^{(2;\Dcal, \Dcal')}_{\alpha, m}
    :=
    \sum^{n_0 -1}_{n = 1}
        \left|
            \Theta^Y
            \binom{s_{m-1}, s_m, s_{m+1}}{u_{n-1}, u_n, u_{n+1}}
        \right|^{\alpha},
\]
for which a calculation in the same vein tells us that there exists $m^* \in \{1, \dotsc, m_0 -1\}$ such that the inequality
\begin{equation}\label{variation inequality alt}
    \Pcal^{(2; \Dcal,\Dcal')}_{\alpha, m^*}
    \leq
    \frac{C(p, \alpha, \theta)}{(m_0 -1)^{\alpha \theta_*}}
    \Vbf^{(2;p)}_{[s,t],[u,v]}(X,Y)^{\alpha}
    ,
\end{equation}
holds, where $C(p, \alpha, \theta_*)$ is the same as in \eqref{variation inequality} and
\[
    \Vbf^{(2;p)}_{[s,t],[u,v]}(X,Y)
    =
    \max_{j,k = 0, \dotsc, N}
        \left\{
            V^{p/(j+1)}_{[s,t]}\left(
                X^j
            \right)
            V^{p/(k+1)}_{[u,v]}\left(
                X^k
            \right)
            V^{q_k, q_j}_{[u,v]\times[s,t]}\left(
                \Rbf^{(2; k, j)}
            \right)
        \right\}
    .
\]
\end{remark}

The proof of this bound draws inspiration from the proof of Lemma 6.4 of \cite{FV11} (original result from \cite{Tow02}) with modifications made to be used in the context of rough paths and jointly controlled paths.

Our strategy now is to obtain the maximal inequality by successively removing well chosen points from the discrete integrals until we are left with a discrete integral over the trivial partition. In the process of removing points, the following intermediate estimates become useful; which we obtain by using the previous bound of Lemma \ref{Theta inequality lemma} and methods from the proof of the same lemma.

\begin{lemma}\label{Delta inequalities lemma}
Suppose that Condition \ref{mixed variation condition} holds. For partitions $\Dcal = \{s_0< \dotsc < s_{m_0} \} \subset [s,t]$ and $\Dcal' = \{u_0 < \dotsc < u_{n_0}\}\subset [u,v]$ define the differences
\begin{alignat*}{2}
    \Delta^{(1; m)}_{\Dcal \times \Dcal'}
    &:=
    \sum_{\Dcal \times \Dcal'}
        \Omega^Y
    -
    \sum_{\Dcal\setminus\{s_{m}\} \times \Dcal'}
        \Omega^Y
    &&=
    \sum_{n=1}^{n_0}
        \Gamma^Y
        \binom{s_{m-1}, s_{m}, s_{m+1}}{u_{n-1}, u_{n}}
    \\[0.5em]
    \Delta^{(2; n)}_{\Dcal \times \Dcal'}
    &:=
    \sum_{\Dcal \times \Dcal'}
        \Omega^Y
    -
    \sum_{\Dcal \times \Dcal'\setminus\{u_n\}}
        \Omega^Y
    &&=
    \sum_{m=1}^{m_0}
        \Gamma^Y
        \binom{s_{m-1}, s_{m}}{u_{n-1}, u_{n}, u_{n+1}}
    .
\end{alignat*}
Then for any $\alpha \in (1/\theta_*, 1)$,
\begin{align}
    \label{Gamma inequality 1}
    \sum_{m = 1}^{m_0 -1}
        \left|
            \Delta^{(1;m)}_{\Dcal \times \{u, v\}}
        \right|^{\alpha}
    &\leq
    \sum^{N}_{j,k = 0}
        \left|
            X^{k+1}_{u, v}
        \right|^{\alpha}
        V^{p/(j+1)}_{[s,t]}
        \left(
            X^{j+1}
        \right)^{\alpha}
        V^{q_j}_{[s,t]}
        \left(
            R^{(2;k,j)}_u
        \right)^{\alpha}
    =:
    \left(
        \eta^{(1;p, \alpha)}_{[s,t],[u,v]}
        (X, Y)
    \right)^{\alpha}
    \\[0.5em]
    \label{Gamma inequality 2}
    \sum_{n = 1}^{n_0 -1}
        \left|
            \Delta^{(2;n)}_{\{s,t\} \times \Dcal'}
        \right|^{\alpha}
    &\leq
    \sum^{N}_{j,k = 0}
        \left|
            X^{j+1}_{s, t}
        \right|^{\alpha}
        V^{p/(k+1)}_{[u,v]}
        \left(
            X^{k+1}
        \right)^{\alpha}
        V^{q_k}_{[u,v]}
        \left(
            R^{(1;j,k)}_s
        \right)^{\alpha}
    =:
    \left(
        \eta^{(2;p, \alpha)}_{[s,t],[u,v]}
        (X, Y)
    \right)^{\alpha}
\end{align}
and there exists constant $C(p,\alpha,\theta_*)$ such that
\begin{align}
    \label{Delta inequality 1}
    \sum^{m_0 -1}_{m = 1}
        \left|
            \Delta^{(1; m)}_{\Dcal \times \Dcal'}
            -
            \Delta^{(1; m)}_{\Dcal \times \{u, v\}}
        \right|^{\alpha}
    &\leq
    C(p,\alpha,\theta_*)
    \zeta\left(
        \alpha \theta_*
    \right)
    \Vbf^{(1;p)}_{[s,t],[u,v]}(X,Y)^{\alpha}
    \\[0.5em]
    \label{Delta inequality 2}
    \sum^{n_0 -1}_{n = 1}
        \left|
            \Delta^{(2; n)}_{\Dcal \times \Dcal'}
            -
            \Delta^{(2; n)}_{\{s,t\} \times \Dcal'}
        \right|^{\alpha}
    &\leq
    C(p,\alpha,\theta_*)
    \zeta\left(
        \alpha \theta_*
    \right)
    \Vbf^{(2;p)}_{[s,t],[u,v]}(X,Y)^{\alpha}
    ,
\end{align}
where $\zeta$ is the Riemann zeta function.
\end{lemma}

\begin{proof}
For \eqref{Gamma inequality 1} and \eqref{Gamma inequality 2} we use similar calculations to those in Lemma \ref{Theta inequality lemma},
\begin{align*}
    \sum_{m = 1}^{m_0 -1}
        \left|
            \Delta^{(1;m)}_{\Dcal \times \{u, v\}}
        \right|^{\alpha}
    &=
    \sum_{m = 1}^{m_0 -1}
        \left|
            \Gamma^Y
            \binom{s_{m-1}, s_{m}, s_{m+1}}{u, v}
        \right|^{\alpha}
    \\[0.5em]
    &\leq
    \sum^{N}_{j,k = 0}
        \sum_{m = 1}^{m_0 -1}
            \left|
                R^{(2;k,j)}_{u; s_{m-1}, s_{m}}
            \right|^{\alpha}
            \left|
                X^{j+1}_{s_{m}, s_{m+1}}
            \right|^{\alpha}
            \left|
                X^{k+1}_{u, v}
            \right|^{\alpha}
    \\[0.5em]
    &\leq
    \sum^{N}_{j,k = 0}
        \left|
            X^{k+1}_{u, v}
        \right|^{\alpha}
        \left(
            \sum_{m = 1}^{m_0 -1}
                \left|
                    X^{j+1}_{s_{m}, s_{m+1}}
                \right|^{p_j}
        \right)^{\alpha/p_j}
        \left(
            \sum_{m = 1}^{m_0 -1}
                \left|
                    R^{(2;k,j)}_{u; s_{m-1}, s_{m}}
                \right|^{q_j}
        \right)^{\alpha/q_j}
    \\[0.5em]
    &\leq
    \sum^{N}_{j,k = 0}
        \left|
            X^{k+1}_{u, v}
        \right|^{\alpha}
        V^{p_j}_{[s,t]}
        \left(
            X^{j+1}
        \right)^{\alpha}
        V^{q_j}_{[s,t]}
        \left(
            R^{(2;k,j)}_u
        \right)^{\alpha}.
\end{align*}
The inequality \eqref{Gamma inequality 2} follows similarly so we omit the calculation.

Likewise the derivation of \eqref{Delta inequality 2} also follows similarly to that of \eqref{Delta inequality 1} so we only show the latter. We first note that for $n \in \{1, \dotsc, n_0 - 1\}$ we have
\[
    \Delta^{(1; m)}_{\Dcal \times \Dcal'}
    -
    \Delta^{(1; m)}_{\Dcal \times \Dcal'\setminus \{u_{n}\}}
    =
    \Theta^Y
    \binom{s_{m-1}, s_{m}, s_{m+1}}{u_{n-1}, u_{n}, u_{n+1}}.
\]

By Lemma \ref{Theta inequality lemma} we may iteratively remove $n^*$ such that \eqref{variation inequality} is satisfied for each new partition formed from these removals. Denote the sequence of $n^*$ by $(n^*_l)_{l=1}^{n_0 - 1}$ and the resulting partitions $\Dcal'_{l} := \Dcal'_{l-1}\setminus \{u_{n^*_l}\}$ with $\Dcal'_0 := \Dcal'$. We obtain a bound on the first term in the above inequality:
\begin{align*}
    \sum^{m_0 -1}_{m = 1}
        \left|
            \Delta^{(1; m)}_{\Dcal \times \Dcal'}
            -
            \Delta^{(1; m)}_{\Dcal \times \{u, v\}}
        \right|^{\alpha}
    &\leq
    \sum^{m_0 -1}_{m = 1}
        \sum^{n_0 - 1}_{l=1}
            \left|
                \Delta^{(1; m)}_{\Dcal \times \Dcal'_{l-1}}
                -
                \Delta^{(1; m)}_{\Dcal \times \Dcal'_{l}}
            \right|^{\alpha}
    \\[0.5em]
    &\leq
    \sum^{n_0 - 1}_{l=1}
        \frac{C(p,\alpha,\theta_*)}{(n_0 -l)^{\alpha \theta_*}}
        \Vbf^{(1;p)}_{[s,t],[u,v]}(X,Y)^{\alpha}
    \\[0.5em]
    &\leq
    C(p,\alpha,\theta_*)
    \zeta\left(
        \alpha \theta_*
    \right)
    \Vbf^{(1;p)}_{[s,t],[u,v]}(X,Y)^{\alpha}
    ,
\end{align*}
where we used Lemma \ref{Theta inequality lemma} in the second line and $\zeta$ is the Riemann zeta function.
\end{proof}

Our first use of these intermediate estimates is to place bounds on comparisons involving discrete integrals with one or both of the one-dimensional partitions in grid-like partition being trivial.

\begin{lemma}\label{endpoints lemma}
Assume that Condition \ref{mixed variation condition} holds. Then for any $\alpha \in (1/\theta_*, 1)$
\begin{align}
    \label{endpoints bound 1}
    \left|
        \sum_{\Dcal \times \{u,v\}}
            \Omega^Y
        -
        \sum_{\{s,t\} \times \{u,v\}}
            \Omega^Y
    \right|
    &\leq
    \zeta\left(
        \frac{1}{\alpha}
    \right)
    \eta^{(1;p, \alpha)}_{[s,t],[u,v]}
        (X, Y)
    \\[0.5em]
    \label{endpoints bound 2}
    \left|
        \sum_{\{s,t\} \times \Dcal'}
            \Omega^Y
        -
        \sum_{\{s,t\} \times \{u,v\}}
            \Omega^Y
    \right|
    &\leq
    \zeta\left(
        \frac{1}{\alpha}
    \right)
    \eta^{(2;p, \alpha)}_{[s,t],[u,v]}
        (X, Y)
\end{align}
and there exists $C'(p, \alpha, \theta_*)$ such that
\begin{equation}\label{endpoints bound}
    \left|
        \sum_{\Dcal \times \Dcal'}
            \Omega^Y
        -
        \sum_{\Dcal \times \{u, v\}}
            \Omega^Y
        -
        \sum_{\{s, t\} \times \Dcal'}
            \Omega^Y
        +
        \sum_{\{s, t\} \times \{u, v\}}
            \Omega^Y
    \right|
    \leq
    C'(p, \alpha, \theta_*)
    \Vbf^{(1;p)}_{[s,t], [u,v]}(X,Y).
\end{equation}
\end{lemma}

\begin{proof}
We prove these inequalities by iteratively removing points and then applying inequalities from Lemma \ref{Delta inequalities lemma}.

Let $\Dcal_0 := \Dcal$ and for $l = 1, \dotsc, m_0 -1$, iteratively define $m_l$ to be the index minimising $\left|\Delta^{(1;m)}_{\Dcal_{l-1} \times \{u,v\}}\right|$ where we recursively define $\Dcal_{l} = \Dcal_{l-1} \setminus \{s_{m_{l}}\}$. By our choice of $m_l$ and applying the inequality \eqref{Gamma inequality 1}, we have
\begin{align*}
    \left|
        \Delta^{(1;m_l)}_{\Dcal_{l-1} \times \{u,v\}}
    \right|^{\alpha}
    &\leq
    \frac{1}{m_0 - l}
    \sum_{s_{m} \in \Dcal_{l-1}}
        \left|
            \Delta^{(1;m)}_{\Dcal_{l-1} \times \{u,v\}}
        \right|^{\alpha}
    \\[0.5em]
    &\leq
    \frac{1}{m_0 - l}
    \left(
        \eta^{(1;p, \alpha)}_{[s,t],[u,v]}
        (X, Y)
    \right)^{\alpha}.
\end{align*}
This then gives us
\begin{align*}
    \left|
        \sum_{\Dcal \times \{u,v\}}
            \Omega^Y
        -
        \sum_{\{s,t\} \times \{u,v\}}
            \Omega^Y
    \right|
    &\leq
    \sum_{l = 1}^{m_0 - 1}
        \left|
            \Delta^{(1;m_{l})}_{\Dcal_{l-1} \times \{u,v\}}
        \right|
    \\[0.5em]
    &\leq
    \zeta\left(
        \frac{1}{\alpha}
    \right)
    \eta^{(1;p, \alpha)}_{[s,t],[u,v]}
    (X, Y).
\end{align*}
The inequality \eqref{endpoints bound 2} follows similarly where we instead successively remove points from $\Dcal'$ to minimise $\left|\Delta^{(2;n)}_{\{s,t\} \times \Dcal_{l-1}'}\right|$, where $\Dcal_{l}'$ are defined analogously.

For \eqref{endpoints bound} we again remove points from $\Dcal$, so define $\Dcal_l$ as before but now we choose $m_l$ to minimise $\left|\Delta^{(1;m)}_{\Dcal_{l-1} \times \Dcal'}- \Delta^{(1;m)}_{\Dcal_{l-1} \times \{u,v\}}\right|$. Similarly to before, we have
\begin{align*}
    \left|
        \Delta^{(1;m_l)}_{\Dcal_{l-1} \times \Dcal'}
        -
        \Delta^{(1;m_l)}_{\Dcal_{l-1} \times \{u,v\}}
    \right|^{\alpha}
    &\leq
    \frac{1}{m_0 - l}
    \sum_{s_{m} \in \Dcal_{l-1}}
        \left|
            \Delta^{(1;m)}_{\Dcal_{l-1} \times \Dcal'}
            -
            \Delta^{(1;m)}_{\Dcal_{l-1} \times \{u,v\}}
        \right|^{\alpha}
    \\[0.5em]
    &\leq
    \frac{C(p,\alpha,\theta_*)}{m_0 - l}
    \zeta\left(
        \alpha \theta_*
    \right)
    \left(
        \Vbf^{(1;p)}_{[s,t],[u,v]}
        (X, Y)
    \right)^{\alpha},
\end{align*}
where we used \eqref{Delta inequality 1} in the second line. We note that
\[
    \sum_{l=0}^{m_0-1}
        \Delta^{(1;m_l)}_{\Dcal_l \times \Dcal'}
        -
        \Delta^{(1;m_l)}_{\Dcal_l \times \{u,v\}}
    =
    \sum_{\Dcal \times \Dcal'}
        \Omega^Y
    -
    \sum_{\Dcal \times \{u, v\}}
        \Omega^Y
    -
    \sum_{\{s, t\} \times \Dcal'}
        \Omega^Y
    +
    \sum_{\{s, t\} \times \{u, v\}}
        \Omega^Y.
\]
Combining these we then arrive at the following inequality:
\begin{align*}
    \left|
        \sum_{l=0}^{m_0-1}
            \Delta^{(1;m_l)}_{\Dcal_l \times \Dcal'}
            -
            \Delta^{(1;m_l)}_{\Dcal_l \times \{u,v\}}
    \right|
    &\leq
    \sum_{l=0}^{m_0-1}
        \left|
            \Delta^{(1;m_l)}_{\Dcal_l \times \Dcal'}
            -
            \Delta^{(1;m_l)}_{\Dcal_l \times \{u,v\}}
        \right|
    \\[0.5em]
    &\leq
    \zeta\left(
        \frac{1}{\alpha}
    \right)
    \left(
        C(p,\alpha,\theta_*)
        \zeta(\alpha\theta_*)
    \right)^{1/\alpha}
    \Vbf^{(1;p)}_{[s,t], [u,v]}(X,Y)
    ,
\end{align*}
which gives us \eqref{Delta inequality 1} with $C'(p, \alpha, \theta_*) = \zeta(1/\alpha)(C(p,\alpha,\theta_*)\zeta(\alpha\theta_*))^{1/\alpha}$.
\end{proof}

\begin{remark}
We may also obtain \eqref{endpoints bound} in a similar fashion by removing points from $\Dcal'$ which would then lead to the inequality
\begin{equation}\label{endpoints bound alt}
    \left|
        \sum_{\Dcal \times \Dcal'}
            \Omega^Y
        -
        \sum_{\Dcal \times \{u, v\}}
            \Omega^Y
        -
        \sum_{\{s, t\} \times \Dcal'}
            \Omega^Y
        +
        \sum_{\{s, t\} \times \{u, v\}}
            \Omega^Y
    \right|
    \leq
    C'(p, \alpha, \theta_*)
    \Vbf^{(2;p)}_{[s,t], [u,v]}(X,Y),
\end{equation}
with the same $C'(p, \alpha, \theta_*)$. This allows us to then write
\begin{equation}\label{endpoints bound final}
    \left|
        \sum_{\Dcal \times \Dcal'}
            \Omega^Y
        -
        \sum_{\Dcal \times \{u, v\}}
            \Omega^Y
        -
        \sum_{\{s, t\} \times \Dcal'}
            \Omega^Y
        +
        \sum_{\{s, t\} \times \{u, v\}}
            \Omega^Y
    \right|
    \leq
    C'(p, \alpha, \theta_*)
    \Vbf^{(p; X, Y)}_{[s,t], [u,v]},
\end{equation}
where
\[
    \Vbf^{(p; X, Y)}_{[s,t], [u,v]}
    =
    \min \left(
        \Vbf^{(1;p)}_{[s,t], [u,v]}(X,Y),
        \Vbf^{(2;p)}_{[s,t], [u,v]}(X,Y)
    \right)
    .
\]
\end{remark}

From these bounds the following maximal inequality is easily obtainable.

\begin{theorem}[Maximal inequality for two-parameter discrete rough integrals]\label{maximal inequality theorem}
Suppose that Condition \ref{mixed variation condition} holds. Then for $\alpha \in (1/\theta_*, 1)$ there exists constant $C''(p, \alpha, \theta_*)$ such that for partitions $\Dcal = \{s_0 < \dotsc < s_{m_0}\} \subset [s,t]$ and $\Dcal' = \{u_0 < \dotsc < u_{n_0}\} \subset [u,v]$ we have the bound
\begin{equation}\label{maximal inequality}
    \left|
        \sum_{\Dcal \times \Dcal'}
            \Omega^Y
        -
        \Omega^Y
        \binom{s,t}{u,v}
    \right|
    \leq
    C''(p, \alpha, \theta_*)
    \left(
        \Vbf^{(p; X, Y)}_{[s,t], [u,v]}
        +
        \eta^{(1;p, \alpha)}_{[s,t],[u,v]}
        (X, Y)
        +
        \eta^{(2;p, \alpha)}_{[s,t],[u,v]}
        (X, Y)
    \right)
\end{equation}
\end{theorem}

\begin{proof}
This follows immediately from Lemma \ref{endpoints lemma} and writing
\begin{alignat*}{2}
    \left|
        \sum_{\Dcal \times \Dcal'}
            \Omega^Y
        -
        \Omega^Y
        \binom{s,t}{u,v}
    \right|
    &\leq &&
    \left|
        \sum_{\Dcal \times \Dcal'}
            \Omega^Y
        -
        \sum_{\Dcal \times \{u, v\}}
            \Omega^Y
        -
        \sum_{\{s, t\} \times \Dcal'}
            \Omega^Y
        +
        \sum_{\{s, t\} \times \{u, v\}}
            \Omega^Y
    \right|
    \\
    & &&+
    \left|
        \sum_{\Dcal \times \{u,v\}}
            \Omega^Y
        -
        \sum_{\{s,t\} \times \{u,v\}}
            \Omega^Y
    \right|
    +
    \left|
        \sum_{\{s,t\} \times \Dcal'}
            \Omega^Y
        -
        \sum_{\{s,t\} \times \{u,v\}}
            \Omega^Y
    \right|.
\end{alignat*}
\end{proof}

To end this section, we give a lemma which provides a bound on the difference of the sum of local approximations when a point is selectively removed from one of the one-dimensional partitions. We note that it is also possible to prove the maximal inequality directly from this lemma.

\begin{lemma}\label{point removal lemma}
Assume Condition \eqref{mixed variation condition} holds. For $\alpha \in (1/\theta_*, 1)$ there exists $m^* \in \{1, \dotsc, m_0 -1 \}$, $n^* \in \{1, \dotsc, n_0 -1 \}$, and constant $C''' = C'''(p, \alpha, \theta_*)$ such that we have
\begin{align}
    \label{removal of a point inequality}
    \left|
        \sum_{\Dcal \times \Dcal'}
            \Omega^Y
        -
        \sum_{\Dcal\setminus\{s_{m^*}\} \times \Dcal'}
            \Omega^Y
    \right|
    &\leq
    C'''
    \left(
        \frac{1}{m_0 - 1}
    \right)^{1/\alpha}
    \left(
        \Vbf^{(1;p)}_{[s,t],[u,v]}(X, Y)
        +
        \eta^{(1;p, \alpha)}_{[s,t],[u,v]}(X, Y)
    \right)
    \\[0.5em]
    \label{removal of a point inequality alt}
    \left|
        \sum_{\Dcal \times \Dcal'}
            \Omega^Y
        -
        \sum_{\Dcal \times \Dcal'\setminus\{u_{n^*}\}}
            \Omega^Y
    \right|
    &\leq
    C'''
    \left(
        \frac{1}{n_0 - 1}
    \right)^{1/\alpha}
    \left(
        \Vbf^{(2;p)}_{[s,t],[u,v]}(X, Y)
        +
        \eta^{(2;p, \alpha)}_{[s,t],[u,v]}(X, Y)
    \right).
\end{align}
\end{lemma}

\begin{proof}
The proof is similar to that of Lemma \ref{endpoints lemma} and again we omit proof of the second inequality which follows similarly to the first. We choose $m^*$ to minimise $\left| \Delta^{(1;m)}_{\Dcal \times \Dcal'}\right|$ and write
\begin{align*}
    \left|
        \Delta^{(1;m^*)}_{\Dcal \times \Dcal'}
    \right|^{\alpha}
    &\leq
    \frac{1}{m_0 - 1}
    \sum_{m=1}^{m_0 - 1}
        \left|
            \Delta^{(1;m)}_{\Dcal \times \Dcal'}
        \right|^{\alpha}
    \\[0.5em]
    &\leq
    \frac{1}{m_0 - 1}
    \sum_{m=1}^{m_0 - 1}
        \left(
            \left|
                \Delta^{(1;m)}_{\Dcal \times \Dcal'}
                -
                \Delta^{(1;m)}_{\Dcal \times \{u,v\}}
            \right|^{\alpha}
            +
            \left|
                \Delta^{(1;m)}_{\Dcal \times \{u,v\}}
            \right|^{\alpha}
        \right)
    \\[0.5em]
    &\leq
    \frac{1}{m_0 - 1}
    \left(
        C(p, \alpha, \theta_*)
        \zeta(\alpha \theta_*)
        \Vbf^{(1;p)}_{[s,t],[u,v]}(X, Y)^{\alpha}
        +
        \eta^{(1;p, \alpha)}_{[s,t],[u,v]}(X, Y)^{\alpha}
    \right)
    \\[0.5em]
    &\leq
    \frac{2^{1-\alpha}}{m_0 - 1}
    \left(
        \left(
            C(p, \alpha, \theta_*)
            \zeta(\alpha \theta_*)
        \right)^{1/\alpha}
        \Vbf^{(1;p)}_{[s,t],[u,v]}(X, Y)
        +
        \eta^{(1;p, \alpha)}_{[s,t],[u,v]}(X, Y)
    \right)^{\alpha},
\end{align*}
where we used Lemma \ref{Delta inequalities lemma} in the second last line and H\"older's inequality on the second. Taking exponent $1/\alpha$ then completes the argument.
\end{proof}

\section{Existence of joint rough integrals and a rough Fubini type theorem}

We now arrive at the main results of this paper, where we show existence and uniqueness of the joint rough integral as well as a rough Fubini type theorem over rectangles. In Section 5 of \cite{GH19} two rough Fubini type theorems are proven (one on the simplex and one on the rectangle) in the $2 \leq p <3$ case under the additional assumption of admitting smooth approximations, which we are able to bypass at the cost of some mild uniformity conditions on the remainders and assumptions of controlled mixed variation.

\begin{theorem}[Existence of rough joint integral]\label{integral existence theorem}
Let $\{Y^{(i;j,k)}\}$ be a jointly $X$-controlled path where $X$ has control $\omega$. For $l = 0, \dotsc, N$, let $p_l:= p/(l+1)$ and assume for $j,k = 0 \dotsc, N$ that $\Rbf^{(1;j,k)}$ have finite $\omega$-controlled $(q_{j},q_{k})$-variation for some $q_l \geq 1$ such that $\frac{1}{p_l} + \frac{1}{q_l} =: \theta_l > 1$.

Furthermore suppose that the following uniformity condition holds: for all $i = 1,2$ and $j,k \in \{0, \dotsc, N\}$, there exists constant $C$ such that the remainders $R^{(i;j,k)}$ satisfy
\begin{equation}\label{uniformity on remainders}
    \sup_{0 \leq u \leq T}
        \left\|
            R^{(i;j,k)}_{u}
        \right\|_{q_k}
    \leq
    C
    .
\end{equation}
Then the integral
\begin{equation}\label{integral existence}
    \int_{[s,t] \times [u,v]}
        Y_{r,r'}
    \;d(X_r, X_{r'})
    :=
    \lim_{|\Dcal \times \Dcal'| \rightarrow 0}
        \sum_{\Dcal \times \Dcal'}
            \Omega^Y
\end{equation}
exists and is such that
\begin{equation}\label{integral bound}
    \left|
        \int_{[s,t] \times [u,v]}
            Y_{r,r'}
        \;d(X_r, X_{r'})
        -
        \Omega^Y
        \binom{s,t}{u,v}
    \right|
    \leq
    C''(p, \alpha, \theta_*)
    \left(
        \Vbf^{(p; X, Y)}_{[s,t], [u,v]}
        +
        \eta^{(p, \alpha)}_{[s,t],[u,v]}
        (X, Y)
    \right),
\end{equation}
where
\[
    \eta^{(p, \alpha)}_{[s,t],[u,v]}
        (X, Y)
    =
    \eta^{(1;p, \alpha)}_{[s,t],[u,v]}
        (X, Y)
    +
    \eta^{(2;p, \alpha)}_{[s,t],[u,v]}
        (X, Y).
\]
\end{theorem}

\begin{proof}
By our regularity assumptions, there exists constant $C$ such that
\begin{align}
    \label{Vbf omega bound}
    \Vbf^{(p;X,Y)}_{[s,t],[u,v]}
    &\leq
    C
    \omega(s,t)^{\theta_*}
    \omega(u,v)^{\theta_*}
    \\[0.5em]
    \label{eta 1 omega bound}
    \eta^{(1;p,\alpha)}_{[s,t],[u,v]}(X,Y)
    &\leq
    C
    \max_k
        \left|
            X^k_{u,v}
        \right|
    \omega(s,t)^{\theta_*}
    \\[0.5em]
    \label{eta 2 omega bound}
    \eta^{(2;p,\alpha)}_{[s,t],[u,v]}(X,Y)
    &\leq
    C
    \max_j
        \left|
            X^j_{s,t}
        \right|
    \omega(u,v)^{\theta_*}
\end{align}
for any $[s,t], [u,v] \subset [0,T]$.

Once existence of \eqref{integral existence} is established, the bound \eqref{integral bound} follows immediately from the maximal inequality of Theorem \ref{maximal inequality theorem}. Since we have a bound from the maximal inequality, we only need to show that given two grid-like partitions $\Dcal \times \Dcal'$ and $\Dcaltilde \times \Dcaltilde'$ that the difference between the sums on the two partitions
\[
    \left|
        \sum_{\Dcal \times \Dcal'}
            \Omega^Y
        -
        \sum_{\Dcaltilde \times \Dcaltilde'}
            \Omega^Y
    \right|
\]
goes to zero as the mesh size $\max \{\Dcal, \Dcal', \Dcaltilde, \Dcaltilde'\}$ goes to zero. Without loss of generality we may assume that $\Dcaltilde \subset \Dcal$ and $\Dcaltilde' \subset \Dcal'$, since otherwise we compare to $\Dcal \cup \Dcaltilde \times \Dcal' \cup \Dcaltilde'$ and use that
\[
    \left|
        \sum_{\Dcal \times \Dcal'}
            \Omega^Y
        -
        \sum_{\Dcaltilde \times \Dcaltilde'}
            \Omega^Y
    \right|
    \leq
    \left|
        \sum_{\Dcal \times \Dcal'}
            \Omega^Y
        -
        \sum_{(\Dcal \cup \Dcaltilde) \times (\Dcal' \cup \Dcaltilde')}
            \Omega^Y
    \right|
    +
    \left|
        \sum_{\Dcaltilde \times \Dcaltilde'}
            \Omega^Y
        -
        \sum_{(\Dcal \cup \Dcaltilde) \times (\Dcal' \cup \Dcaltilde')}
            \Omega^Y.
    \right|
\]
Under this assumption of $\Dcal \times \Dcal'$ refining $\Dcaltilde \times \Dcaltilde'$, we start by splitting up this difference similarly to what we did in Theorem \ref{maximal inequality theorem},
\begin{alignat*}{2}
    \left|
        \sum_{\Dcal \times \Dcal'}
            \Omega^Y
        -
        \sum_{\Dcaltilde \times \Dcaltilde'}
            \Omega^Y
    \right|
    &\leq &&
    \left|
        \sum_{\Dcal \times \Dcal'}
            \Omega^Y
        -
        \sum_{\Dcal \times \Dcaltilde'}
            \Omega^Y
        -
        \sum_{\Dcaltilde \times \Dcal'}
            \Omega^Y
        +
        \sum_{\Dcaltilde \times \Dcaltilde'}
            \Omega^Y
    \right|
    \\
    & &&+
    \left|
        \sum_{\Dcal \times \Dcaltilde'}
            \Omega^Y
        -
        \sum_{\Dcaltilde \times \Dcaltilde'}
            \Omega^Y
    \right|
    +
    \left|
        \sum_{\Dcaltilde \times \Dcal'}
            \Omega^Y
        -
        \sum_{\Dcaltilde \times \Dcaltilde'}
            \Omega^Y
    \right|
    .
\end{alignat*}
Label each of these three terms as $I_1, I_2, I_3$ respectively. As a control $\omega$ is continuous and thus uniformly continuous on $[0,T] \times [0,T]$. So for $\epsilon > 0$, let $\delta_{\epsilon} = \sup\{ \omega(a,b) : a,b \in [0,T], |b-a| < \epsilon \}$ and suppose that $\max(|\Dcal|, |\Dcal'|) < \epsilon$.

Writing $\Dcaltilde = \{ s_0 < \dotsc < s_{m_0}\}$ and $\Dcaltilde' = \{ u_0 < \dotsc u_{n_0}\}$, we may view $\Dcal$ and $\Dcal'$ as collections of partitions $(\Dcal_m)_{m=1}^{m_0}$ and $(\Dcal_n')_{n=1}^{n_0}$ where $\Dcal_m = \Dcal \cap [s_{m-1}, s_m]$ and $\Dcal_n' = \Dcal' \cap [u_{n-1}, u_n]$. Under this representation, we rewrite $I_1$ as
\begin{align*}
    I_1
    &=
    \left|
        \sum_{m=1}^{m_0}
        \sum_{n=1}^{n_0}
            \left(
                \sum_{\Dcal_m \times \Dcal_n'}
                    \Omega^Y
                -
                \sum_{\Dcal_m \times \{u_{n-1}, u_n\}}
                    \Omega^Y
                -
                \sum_{\{s_{m-1}, s_m\} \times \Dcal_n'}
                    \Omega^Y
                +
                \sum_{\{s_{m-1}, s_m\} \times \{u_{n-1}, u_n\}}
                    \Omega^Y
            \right)
    \right|
    \\[0.5em]
    &\leq
    \sum_{m=1}^{m_0}
    \sum_{n=1}^{n_0}
        \left|
            \sum_{\Dcal_m \times \Dcal_n'}
                \Omega^Y
            -
            \sum_{\Dcal_m \times \{u_{n-1}, u_n\}}
                \Omega^Y
            -
            \sum_{\{s_{m-1}, s_m\} \times \Dcal_n'}
                \Omega^Y
            +
            \sum_{\{s_{m-1}, s_m\} \times \{u_{n-1}, u_n\}}
                \Omega^Y
        \right|
    \\[0.5em]
    &\leq
    \sum_{m=1}^{m_0}
    \sum_{n=1}^{n_0}
        C'(p, \alpha, \theta_*)
        \Vbf^{(p; X, Y)}_{[s_{m-1}, s_m], [u_{n-1}, u_n]}
    \\[0.5em]
    &\leq
    C'(p, \alpha, \theta_*)
    \sum_{m=1}^{m_0}
    \sum_{n=1}^{n_0}
        C
        \omega(s_{m-1},s_m)^{\theta_*}
        \omega(u_{n-1},u_n)^{\theta_*}
    \\[0.5em]
    &=
    O \left(
        |\delta_{\epsilon}|^{2(\theta_* - 1)}|
    \right),
\end{align*}
where we used Lemma \ref{endpoints bound} on the third line, the inequality \eqref{Vbf omega bound} on the penultimate line.

For the term $I_2$ we use the same division into subpartitions to obtain the following:
\begin{align*}
    I_2
    &\leq
    \sum_{m=1}^{m_0}
        \left|
            \sum_{\Dcal_m \times \Dcaltilde}
                \Omega^Y
            -
            \sum_{\{s_{m-1}, s_m\} \times \Dcaltilde}
                \Omega^Y
        \right|
    \\[0.5em]
    &\leq
    \sum_{m=1}^{m_0}
        C'''
        \zeta\left(
            \frac{1}{\alpha}
        \right)
        \left(
            \Vbf^{(p; X, Y)}_{[s_{m-1}, s_m], [u,v]}
            +
            \eta^{(1;p,\alpha)}_{[s_{m-1},s_m],[u,v]}(X,Y)
        \right)
    \\[0.5em]
    &=
    O \left(
        |\delta_{\epsilon}|^{\theta_* - 1}
    \right),
\end{align*}
where in the second line we used Lemma \ref{point removal lemma} and successively removed points from $\Dcal_m$ to reduce it to $\{s_{m-1}, s_m\}$. A similar calculation on $I_3$ then yields
\begin{align*}
    I_3
    &\leq
    \sum^{n_0}_{n=1}
        \left|
            \sum_{\Dcal \times \Dcal_n}
                \Omega^Y
            -
            \sum_{\Dcal \times \{u_{n-1}, u_n\}}
                \Omega^Y
        \right|
    \\[0.5em]
    &\leq
    O \left(
        |\delta_{\epsilon}|^{\theta_* - 1}
    \right).
\end{align*}
Since $\delta_{\epsilon} \rightarrow 0$ as $\epsilon \rightarrow 0$, we deduce that the limit \eqref{integral existence} is well defined, completing the proof.
\end{proof}

Following similar methods to those used in proving existence of the joint integral, we are also able to give a rough Fubini type theorem.

\begin{theorem}[Rough Fubini type theorem on the rectangle]\label{rough Fubini theorem}
Suppose that all the conditions in Theorem \ref{integral existence theorem} hold. Then the double integrals are equal to the joint integral \eqref{integral existence},
\begin{equation}\label{rough Fubini}
    \int_s^t
        \int_u^v
            Y_{r,r'}
        \;dX_{r'}
    \;dX_{r})
    =
    \int_{[s,t] \times [u,v]}
        Y_{r,r'}
    \;d(X_r, X_{r'})
    =
    \int_u^v
        \int_s^t
            Y_{r,r'}
        \;dX_{r}
    \;dX_{r'}
    .
\end{equation}
\end{theorem}

\begin{proof}
We prove this by showing that the local approximations $\sum_{\Dcal \times \Dcal'} \Omega^Y$ converge to the double integrals as $|\Dcal \times \Dcal'| = \max(|\Dcal|, |\Dcal'|) \rightarrow 0$. Recall that $Z_s^{(j)} = \int_u^v Y_{s, r}^{(1;j,0)} \; dX_r$ are $X$-controlled paths and by definition
\begin{align*}
    Z^{(j)}_s
    \Big(
        X^{j+1}_{s, t}
    \Big)
    &=
    \lim_{|\Dcal'| \rightarrow 0}
        \sum_{\{u_{n-1}, u_{n}\} \subset \Dcal'}
            \left(
                \sum_{k = 0}^{N}
                    Y_{s, u_{n-1}}^{(1;j,k)}
                    \Big(
                        X^{k+1}_{u_{n-1}, u_{n}}
                    \Big)
                    \Big(
                        X^{j+1}_{s, t}
                    \Big)
            \right)
    \\[0.5em]
    \Ztilde^{(j)}_u
    \Big(
        X^{j+1}_{u, v}
    \Big)
    &=
    \lim_{|\Dcal| \rightarrow 0}
        \sum_{\{s_{m-1}, s_{m}\} \subset \Dcal}
            \left(
                \sum_{k = 0}^{N}
                    Y_{s_{m-1}, u}^{(2;j,k)}
                    \Big(
                        X^{k+1}_{s_{m-1}, s_{m}}
                    \Big)
                    \Big(
                        X^{j+1}_{u, v}
                    \Big)
            \right)
\end{align*}
which gives us
\begin{align*}
    \sum_{j=0}^N
        Z^{(j)}_s
        \Big(
            X^{j+1}_{s, t}
        \Big)
    =
    \lim_{|\Dcal'| \rightarrow 0}
        \sum_{\{s,t\} \times \Dcal'}
            \Omega^Y,
    \qquad
    \sum_{j=0}^N
        \Ztilde^{(j)}_u
        \Big(
            X^{j+1}_{u, v}
        \Big)
    =
    \lim_{|\Dcal| \rightarrow 0}
        \sum_{\Dcal \times \{u,v\}}
            \Omega^Y.
\end{align*}
Let $\Dcal = \{s_0 < \dotsc < s_{m_0}\} \subset [s,t]$ and $\Dcal' = \{u_0 < \dotsc < u_{n_0}\} \subset [u,v]$.
\begin{alignat*}{2}
    \left|
        \int_s^t
            \int_u^v
                Y_{r,r'}
            \;dX_{r'}
        \;dX_{r}
        -
        \sum_{\Dcal \times \Dcal'}
            \Omega^Y
    \right|
    &\leq &&
    \left|
        \int_s^t
            \int_u^v
                Y_{r,r'}
            \;dX_{r'}
        \;dX_{r}
        -
        \sum_{m=1}^{m_0}
            \sum_{j = 0}^{N}
                Z^{(j)}_{s_{m-1}}
                \left(
                    X^{j+1}_{s_{m-1}, s_{m}}
                \right)
    \right|
    \\[0.5em]
    &
    &&+
    \left|
        \sum_{m=1}^{m_0}
            \sum_{j = 0}^{N}
                Z^{(j)}_{s_{m-1}}
                \left(
                    X^{j+1}_{s_{m-1}, s_{m}}
                \right)
        -
        \sum_{\Dcal \times \Dcal'}
            \Omega^Y
    \right|
    \\[1em]
    \left|
        \int_u^v
            \int_s^t
                Y_{r,r'}
            \;dX_{r}
        \;dX_{r'}
        -
        \sum_{\Dcal \times \Dcal'}
            \Omega^Y
    \right|
    &\leq &&
    \left|
        \int_u^v
            \int_s^t
                Y_{r,r'}
            \;dX_{r}
        \;dX_{r'}
        -
        \sum_{n=1}^{n_0}
            \sum_{j = 0}^{N}
                \Ztilde^{(j)}_{u_{n-1}}
                \left(
                    X^{j+1}_{u_{n-1}, s_{n}}
                \right)
    \right|
    \\[0.5em]
    &
    &&+
    \left|
        \sum_{n=1}^{n_0}
            \sum_{j = 0}^{N}
                \Ztilde^{(j)}_{u_{n-1}}
                \left(
                    X^{j+1}_{u_{n-1}, u_{n}}
                \right)
        -
        \sum_{\Dcal \times \Dcal'}
            \Omega^Y
    \right|
\end{alignat*}
From the definition of controlled path integration, the first term in each of the inequalities tends to zero as the mesh size $|\Dcal|$ ($|\Dcal'|$ respectively) goes to zero.

Using the same notation as in Theorem \ref{integral existence theorem}, for $\epsilon > 0$ and $|\Dcal \times \Dcal'| < \epsilon$, the bounds on $I_3$ and $I_2$ from Theorem \ref{integral existence theorem} gives us existence of a constants $C(s,t) = C(p, \alpha, X, Y, s, t)$ and $C(u,v) = C(p, \alpha, X, Y, u, v)$ such that
\begin{align*}
    \left|
        \sum_{m=1}^{m_0}
            \sum_{j = 0}^{N}
                Z^{(j)}_{s_{m-1}}
                \left(
                    X^{j+1}_{s_{m-1}, s_{m}}
                \right)
        -
        \sum_{\Dcal \times \Dcal'}
            \Omega^Y
    \right|
    &=
    \lim_{|\Dcaltilde'| \rightarrow 0}
        \left|
            \sum_{\Dcal \times \Dcaltilde'}
                \Omega^Y
            -
            \sum_{\Dcal \times \Dcal'}
                \Omega^Y
        \right|
    \leq
    C(s,t)
    \;|\delta_{\epsilon}|^{\theta_*-1}
    \\[0.5em]
    \left|
        \sum_{n=1}^{n_0}
            \sum_{j = 0}^{N}
                \Ztilde^{(j)}_{u_{n-1}}
                \left(
                    X^{j+1}_{u_{n-1}, u_{n}}
                \right)
        -
        \sum_{\Dcal \times \Dcal'}
            \Omega^Y
    \right|
    &=
    \lim_{|\Dcaltilde| \rightarrow 0}
        \left|
            \sum_{\Dcaltilde \times \Dcal'}
                \Omega^Y
            -
            \sum_{\Dcal \times \Dcal'}
                \Omega^Y
        \right|
    \leq
    C(u,v)
    \;|\delta_{\epsilon}|^{\theta_*-1}
    .
\end{align*}
With this we have established that our local approximation $\sum_{\Dcal \times \Dcal'} \Omega^Y$ converges to both iterated integrals, which gives us \eqref{rough Fubini} as required.
\end{proof}

The rough Fubini type theorems derived in \cite{GH19} cover the case where the controlling rough paths have finite $p$-variation for $2 \leq p < 3$, and are used as tools in order to prove H\"ormander's theorem for a class of SPDEs. In comparison with the Fubini type theorems of \cite{GH19}, we do not cover integrals over the simplex, but are able to generalise integrals over rectangles to the case where the controlling rough path is geometric with arbitrary $p$-variation. By introducing the third joint integral, the proof given here draws more parallels with the classical Fubini's theorem, and allows us to bypass the technical condition of smooth approximability.

Next we move on to a stability result on the double integral, analogous to the stability result of the one parameter rough integral in the case $2 \leq p < 3$ as presented in Theorem 4.17 of \cite{FH14}. This result in essence tells us that if two controlled paths are close in some sense, then the double integrals are also close. In order to show this we rely on the following type of estimate:

\begin{lemma}\label{difference of controlled paths lemma}
For two $\omega$-controlled geometric $p$-rough paths $X$ and $\Xtilde$, let $(A^0, \dotsc, A^N)$ be an $X$-controlled path on a Banach space $E$ with remainders $\{R^j\}_{j=0}^N$ and let $(\Atilde^0, \dotsc, \Atilde^N)$ be an $\Xtilde$-controlled path with remainders $\{ \Rtilde^j \}_{j=0}^N$. Define
\[
    \Xi^{(A, \Atilde)}_{s,t}
    =
    \sum^{N}_{j=0}
        A^{j}_s
        \big(
            X^{j+1}_{s,t}
        \big)
        -
        \Atilde^{j}_s
        \big(
            \Xtilde^{j+1}_{s,t}
        \big).
\]
Let $p_j = p/(j+1)$ and $q_j > 1$ such that $\frac{1}{p_j} + \frac{1}{q_j} = \theta_j > 1$, and the remainders $R^j, \Rtilde^j$ have finite $\omega$-controlled $q_j$-variation. Writing $\theta_* = \min_{j} \theta_j$, and $\| \cdot \|_{\alpha} = \| \cdot \|_{\alpha, \omega}$, we have that
\begin{equation} \label{difference of controlled paths}
    \left\|
        \delta \Xi^{(A, \Atilde)}
    \right\|_{1/\theta_*}
    \leq
    \sum_{j=0}^{N}
        C(\theta_j, T)
        \left(
            \|
                R^j
            \|_{q_j}
            \|
                X^{j+1}
                -
                \Xtilde^{j+1}
            \|_{p_j}
            +
            \|
                R^j - \Rtilde^j
            \|_{q_j}
            \|
                \Xtilde^{j+1}
            \|_{p_j}
        \right)
        ,
\end{equation}
where $C(\theta_j, T) = \omega(0,T)^{\theta_j - \theta_*}$.
\end{lemma}

\begin{proof}
Using basic estimates of the form $|a_1 b_1 - a_2 b_2| \leq |a_1| |b_1 - b_2| + |a_1 - a_2| |b_2|$ in conjuncture with Lemma \ref{controlled path identity lemma} we have
\begin{align*}
    \left|
        \delta \Xi^{(A, \Atilde)}_{s,s',t}
    \right|
    &=
    \left|
        \sum^N_{j=0}
            R^{j}_{s,s'}
            \big(
                X^{j+1}_{s',t}
            \big)
            \Rtilde^{j}_{s,s'}
            \big(
                \Xtilde^{j+1}_{s',t}
            \big)
    \right|
    \\[0.5em]
    &\leq
    \sum^N_{j=0}
        \left(
            \left\|
                R^{j}
            \right\|_{q_j}
            \left\|
                X^{j+1}
                -
                \Xtilde^{j+1}
            \right\|_{p_j}
            +
            \left\|
                R^{j}
                -
                \Rtilde^{j}
            \right\|_{q_j}
            \left\|
                \Xtilde^{j+1}
            \right\|_{p_j}
        \right)
        \omega(s,s')^{1/q_j}
        \omega(s',t)^{1/p_j}
    \\[0.5em]
    &\leq
    \sum^N_{j=0}
        \left(
            \left\|
                R^{j}
            \right\|_{q_j}
            \left\|
                X^{j+1}
                -
                \Xtilde^{j+1}
            \right\|_{p_j}
            +
            \left\|
                R^{j}
                -
                \Rtilde^{j}
            \right\|_{q_j}
            \left\|
                \Xtilde^{j+1}
            \right\|_{p_j}
        \right)
        \omega(s,t)^{\theta_j}
    ,
\end{align*}
from which \eqref{difference of controlled paths} becomes apparent.
\end{proof}

\begin{theorem}[Stability of double rough integrals]
\label{stability theorem}
Let $X, \Xtilde$ be two geometric $p$-rough paths on $V$ with control $\omega$. Consider two jointly controlled paths $\{Y^{(i;j,k)}\} \in \Dscr_X^p$ and $\{\Ytilde^{(i;j,k)}\} \in \Dscr_{\Xtilde}^p$ with first order remainders $\{R^{(i;j,k)}\}, \{\Rtilde^{(i;j,k)}\}$ and second order remainders $\{\Rbf^{(i;j,k)}\}, \{\Rbftilde^{(i;j,k)}\}$ respectively. Fix some $[s,t], [u,v] \subset [0,T]$. Suppose that these two paths satisfy the conditions of Theorem \ref{integral existence theorem} for some $q_j$, and define the "distance" between two jointly controlled paths
\begin{alignat*}{2}
    d_{X, \Xtilde}(Y, \Ytilde)
    &:=
    \sum^{N}_{j=0}
    \sum^{N}_{k=0}
    &&
        \left|
            Y^{(1;j,k)}_{s,u}
            -
            \Ytilde^{(1;j,k)}_{s,u}
        \right|
        +
        \left\|
            R^{(1;j,k)}_{s}
            -
            \Rtilde^{(1;j,k)}_{s}
        \right\|_{q_k}
        +
        \left\|
            R^{(2;k,j)}_{u}
            -
            \Rtilde^{(2;k,j)}_{u}
        \right\|_{q_j}
    \\
    &
    &&
        +
        \left\|
            \Rbf^{(1;j,k)}
            -
            \Rbftilde^{(1;j,k)}
        \right\|_{q_j, q_k}
    .
\end{alignat*}
The double integrals are such that
\begin{equation}
    \label{stability bound}
    \left|
        \int_{[s,t] \times [u,v]}
            Y_{r,r'}
        \;d(X_r, X_{r'})
        -
        \int_{[s,t] \times [u,v]}
            \Ytilde_{r,r'}
        \;d(\Xtilde_r, \Xtilde_{r'})
    \right|
    \leq
    C
    \left(
        \|X - \Xtilde\|_{p\variation}
        +
        d_{X, \Xtilde}(Y, \Ytilde)
    \right)
    ,
\end{equation}
where $C$ depends on $X, \Xtilde, Y, \Ytilde$, and remainders, at initial and terminal times $s,t,u,v$.
\end{theorem}

\begin{remark}
As previously mentioned, this result is an adaptation of the one-parameter stability result for rough integrals as presented in Theorem 4.17 \cite{FH14}, which similarly defines a "distance" between two controlled paths. We thus make the same remark to note that while we refer to this as a distance, this is not a true metric as the two paths will, in general, live in two different spaces.
\end{remark}

\begin{proof}
We prove this by making estimates on the double integrals and using the rough Fubini theorem established earlier. Define the rough integrals
\begin{equation*}
    Z^{(j)}_r
    =
    \int_u^v
        Y^{(1;j,0)}_{r,r'}
    \;dX_{r'}
    ,
    \qquad
    \Ztilde^{(j)}_r
    =
    \int_u^v
        \Ytilde^{(1;j,0)}_{r,r'}
    \;d\Xtilde_{r'}
    .
\end{equation*}
Recall that by Lemma \ref{integrals are controlled paths} that $(Z^{(j)})$ and $(\Ztilde^{(j)})$ are $X$-controlled and $\Xtilde$-controlled paths respectively with remainders
\begin{align*}
    R^{(Z;j)}_{s,t}
    =
    \int^v_u
        R^{(1;j,0)}_{s,t;r'}
    \;dX_{r'}
    ,
    \qquad
    \Rtilde^{(\Ztilde;j)}_{s,t}
    =
    \int^v_u
        \Rtilde^{(1;j,0)}_{s,t;r'}
    \;dX_{r'}
    ,
\end{align*}
respectively. Define the double integrals
\begin{equation*}
    I
    =
    \int_s^t
        Z^{(0)}_{r}
    \;dX_{r}
    ,
    \qquad
    \Itilde
    =
    \int_s^t
        \Ztilde^{(0)}_{r}
    \;dX_{r}.
\end{equation*}
and the quantity
\[
    \Xi^{Z, \Ztilde}_{s,t}
    =
    \sum^N_{j=0}
        Z^{(j)}_s
        \Big(
            X^{j+1}_{s,t}
        \Big)
        -
        \Ztilde^{(j)}_s
        \Big(
            \Xtilde^{j+1}_{s,t}
        \Big)
    .
\]
Applying Lemma \ref{sewing lemma} where we take $\beta = \theta_*$ and $\Xi = \Xi^{Z, \Ztilde}$, we have
\begin{align*}
    \big|
        I - \Itilde
    \big|
    \leq
    \big|
        \Xi^{Z,\Ztilde}_{s,t}
    \big|
    +
    \zeta(\theta_*)
    \;
    \omega(s,t)^{\theta_*}
    \big\|
        \delta \Xi^{Z, \Ztilde}
    \big\|_{1/\theta_*}
    .
\end{align*}
Using elementary estimates of the form $|a_1 b_1 - a_2 b_2| \leq |a_1||b_1 - b_2| + |a_1 - a_2||b_2|$, for the first term we have
\begin{align*}
    \big|
        \Xi^{Z,\Ztilde}_{s,t}
    \big|
    &\leq
    \sum_{j=0}^N
        \left|
            Z^{(j)}_s
        \right|
        \big|
            X^{j+1}_{s,t}
            -
            \Xtilde^{j+1}_{s,t}
        \big|
        +
        \left|
            Z^{(j)}_s
            -
            \Ztilde^{(j)}_s
        \right|
        \big|
            \Xtilde^{j+1}_{s,t}
        \big|.
\end{align*}
By the sewing lemma again, we have the bound
\[
    \left|
        Z^{(j)}_s
        -
        \Ztilde^{(j)}_s
    \right|
    \leq
    \left|
        \Delta^{(j)}_{s;u,v}
    \right|
    +
    \zeta(\theta_*)
    \big\|
        \delta \Delta^{(j)}_s
    \big\|_{1/\theta_*}
    ,
\]
where
\[
    \Delta^{(j)}_{s;u,v}
    =
    \sum^N_{k=0}
        Y^{(1;j,k)}_{s,u}
        \big(
            X^{k+1}_{u,v}
        \big)
        -
        \Ytilde^{(1;j,k)}_{s,u}
        \big(
            \Xtilde^{k+1}_{u,v}
        \big)
    .
\]
Similarly to before we now have
\[
    \left|
        \Delta^{(j)}_{s;u,v}
    \right|
    \leq
    \sum^N_{k=0}
        \left|
            Y^{(1;j,k)}_{s,u}
        \right|
        \big|
            X^{k+1}_{u,v}
            -
            \Xtilde^{k+1}_{u,v}
        \big|
        +
        \left|
            Y^{(1;j,k)}_{s,u}
            -
            \Ytilde^{(1;j,k)}_{s,u}
        \right|
        \big|
            \Xtilde^{k+1}_{u,v}
        \big|
\]
and using Lemma \ref{difference of controlled paths lemma} we know that
\begin{align*}
    \left\|
        \delta \Delta^{(j)}_s
    \right\|_{1/\theta_*}
    &\leq
    \sum_{k=0}^{N}
        C(\theta_k, T)
        \left(
            \|
                R^{(1;j,k)}_s
            \|_{q_k}
            \|
                X^{k+1}
                -
                \Xtilde^{k+1}
            \|_{p_k}
            +
            \|
                R^{(1;j,k)}_s - \Rtilde^{(1;j,k)}_s
            \|_{q_k}
            \|
                \Xtilde^{k+1}
            \|_{p_k}
        \right)
\end{align*}
which then gives us that
\[
    \big|
        \Xi^{Z,\Ztilde}_{s,t}
    \big|
    \lesssim
    \|
        X - \Xtilde
    \|_{p\variation}
    +
    d_{X, \Xtilde} (Y, \Ytilde)
    .
\]
A straightforward application of Lemma \ref{difference of controlled paths lemma} yields that
\[
    \big\|
        \delta \Xi^{Z, \Ztilde}
    \big\|_{1/\theta_*}
    \leq
    \sum_{j=0}^{N}
        C(\theta_j, T)
        \left(
            \|
                R^{(Z;j)}
            \|_{q_j}
            \|
                X^{j+1}
                -
                \Xtilde^{j+1}
            \|_{p_j}
            +
            \|
                R^{(Z;j)} - \Rtilde^{(\Ztilde;j)}
            \|_{q_j}
            \|
                \Xtilde^{j+1}
            \|_{p_j}
        \right)
    ,
\]
so we now find bounds on the integral remainders $R^{(Z;j)}$ and $\Rtilde^{(\Ztilde;j)}$. These remainders are themselves rough integrals and so using the sewing lemma again we have
\[
    \left|
        R^{(Z;j)}_{s,t}
        -
        \Rtilde^{(\Ztilde;j)}_{s,t}
    \right|
    \leq
    \left|
        \Xi^{(R, \Rtilde; j)}_{s,t;u,v}
    \right|
    +
    \zeta(\theta_*)
    \left\|
        \delta \Xi^{(R, \Rtilde; j)}_{s,t}
    \right\|_{1/\theta_*}
\]
where $\Xi^{(R, \Rtilde; j)}_{s,t;u,v}$ are defined by
\[
    \Xi^{(R, \Rtilde; j)}_{s,t;u,v}
    =
    \sum_{k=0}^{N}
        R^{(1;j,k)}_{s,t;u}
        \big(
            X^{k+1}_{u,v}
        \big)
        -
        \Rtilde^{(1;j,k)}_{s,t;u}
        \big(
            \Xtilde^{k+1}_{u,v}
        \big)
    .
\]
We now repeat the same arguments as earlier to give us
\begin{align*}
    \left|
        \Xi^{(R, \Rtilde; j)}_{s,t;u,v}
    \right|
    &\leq
    \sum_{k=0}^N
        \left|
            R^{(2;k,j)}_{u;s,t}
        \right|
        \big|
            X^{k+1}_{u,v}
            -
            \Xtilde^{k+1}_{u,v}
        \big|
        +
        \left|
            R^{(2;k,j)}_{u;s,t}
            -
            \Rtilde^{(2;k,j)}_{u;s,t}
        \right|
        \big|
            \Xtilde^{k+1}_{u,v}
        \big|
    ,
    \\[0.5em]
    \left\|
        \delta \Xi^{(R, \Rtilde; j)}_{s,t}
    \right\|_{1/\theta_*}
    &\leq
    \sum_{k=0}^{N}
    C_k
        \left(
            \|
                \Rbf^{(1;j,k)}_{s,t}
            \|_{q_k}
            \|
                X^{k+1}
                -
                \Xtilde^{k+1}
            \|_{p_k}
            +
            \|
                \Rbf^{(1;j,k)}_{s,t}
                -
                \Rbftilde^{(1;j,k)}_{s,t}
            \|_{q_k}
            \|
                \Xtilde^{k+1}
            \|_{p_k}
        \right)
\end{align*}
for constants $C_k = C(\theta_k, T)$, which we then use to deduce that
\begin{align*}
    \left\|
        R^{(Z;j)}
        -
        \Rtilde^{(\Ztilde;j)}
    \right\|_{q_j}
    &\leq
    \sum_{k=0}^N
        \big\|
            R^{(2;k,j)}_u
        \big\|_{q_j}
        \big|
            X^{k+1}_{u,v}
            -
            \Xtilde^{k+1}_{u,v}
        \big|
        +
        \big\|
            R^{(2;k,j)}_u
            -
            \Rtilde^{(2;k,j)}_u
        \big\|_{q_j}
        \big|
            \Xtilde^{k+1}_{u,v}
        \big|
    \\
    &\quad
    +
    \sum_{k=0}^{N}
    C_k
        \|
            \Rbf^{(1;j,k)}
        \|_{q_j, q_k}
        \|
            X^{k+1}
            -
            \Xtilde^{k+1}
        \|_{p_k}
    \\
    &\quad
    +
    \sum_{k=0}^{N}
    C_k
        \|
            \Rbf^{(1;j,k)}
            -
            \Rbftilde^{(1;j,k)}
        \|_{q_j, q_k}
        \|
            \Xtilde^{k+1}
        \|_{p_k}
    \\[0.5em]
    &\lesssim
    \|
        X - \Xtilde
    \|_{p\variation}
    +
    d_{X, \Xtilde} (Y, \Ytilde)
    .
\end{align*}
Combining these components we then arrive at \eqref{stability bound} as required.
\end{proof}

\section{Paths controlled by two different rough paths}

In the earlier sections we investigated two-parameter rough integration in the case where we have a single controlling rough path. Here we now take the same ideas and generalise to two-parameter rough integration with respect to two potentially different controlling rough paths. As before, we first begin by defining a class of paths that are controlled by two geometric rough paths $X$ and $\Xtilde$ in order to integrate against them.

\begin{definition}
Let $X = (1, X^1, X^2, \dotsc, X^{\pfloor})$ be an $\omega$-controlled geometric $p$-rough path on a Banach space $V$ and $\Xtilde = (1, \Xtilde^1, \Xtilde^2, \dotsc, \Xtilde^{\ptilfloor})$ an $\omegatil$-controlled geometric $\ptilde$-rough path on $\Vtilde$. Let $N = \pfloor - 1$ and $\Ntilde = \ptilfloor - 1$.

A two parameter path $Y: [0,T]^2 \rightarrow E$ on a Banach space $E$ is jointly $(X, \Xtilde)$-controlled if the following conditions are satisfied:
\begin{enumerate}[topsep=0pt]
    \item Let $Y_{s, \cdot}^{(1;0,0)} = Y_{s, \cdot}$ and $Y_{\cdot, u}^{(2;0,0)} = Y_{\cdot, u}$. For $j= 0, \dotsc, N$, $k = 0, \dotsc, \Ntilde$, and every $s, u \in [0,T]$, there exists
    \begin{align*}
        Y^{(1; j, k)}_{s, \cdot}
        &:
        [0, T]
        \rightarrow
        \Hom\left(
            \Vtilde^{\otimes k},
            \Hom(V^{\otimes j}, E)
        \right)
        \\
        Y^{(2; k, j)}_{\cdot, u}
        &:
        [0, T]
        \rightarrow
        \Hom\left(
            V^{\otimes j},
            \Hom(\Vtilde^{\otimes k}, E)
        \right)
    \end{align*}
    such that $\left( Y^{(1;j,0)}_{s, \cdot}, \dotsc, Y^{(1;j, N)}_{s, \cdot} \right)$ is an $\Xtilde$-controlled path and $\left( Y^{(2;k,0)}_{\cdot, u}, \dotsc, Y^{(2;k, \Ntilde)}_{\cdot, u} \right)$ is an $X$-controlled path.
    
    \item For all $j= 0, \dotsc, N$, $k = 0, \dotsc, \Ntilde$, and $s, u \in [0,T]$, the derivatives $Y^{(1;j,k)}$ and $Y^{(2;k,j)}$ satisfy the symmetry condition
    \begin{equation}\label{derivative symmetry 2}
        Y^{(1;j,k)}_{s,u}(y)(x)
        =
        Y^{(2;k,j)}_{s,u}(x)(y)
    \end{equation}
    for $x \in V^{\otimes j}$ and $y \in \Vtilde^{\otimes k}$.
\end{enumerate}
The collection $\left\{Y^{(i; j, k)} \right\}$ defines the jointly $(X, \Xtilde)$-controlled path, and we denote this class of paths by $\Dscr_{X, \Xtilde}^{p, \ptilde} \left( [0,T]^2; E \right)$.
\end{definition}

Typically we will take $p = \ptilde$, $V = \Vtilde$, and $\omega = \omegatil$, although it is not necessary. As before we also have remainders $R^{(i;j,k)}$, which instead are such that
\begin{align*}
    Y_{s, v}^{(1;j,k)}
    &=
    \sum^{\Ntilde-k}_{l = 0}
        Y_{s, u}^{(1;j,k+l)}
        \Big(
            \Xtilde^l_{u,v}
        \Big)
    +
    R^{(1;j,k)}_{s;u,v}, \\[0.5em]
    Y_{t, u}^{(2;k,j)}
    &=
    \sum^{N-j}_{m = 0}
        Y_{s, u}^{(2;k, j+m)}
        \Big(
            X^m_{s,t}
        \Big)
    +
    R^{(2;k,j)}_{u;s,t},
\end{align*}
and which now satisfy
\begin{align*}
    \sup_{0 \leq s < t \leq T}
        \left|
            \frac{
                R^{(1;j,k)}_{s;u,v}
            }{
                \omegatil(u,v)^{
                    \left(
                        \ptilfloor - k
                    \right)
                    /\ptilde
                }
            }
        \right|
    &< \infty
    ,
    \qquad
    \sup_{0 \leq s < t \leq T}
        \left|
            \frac{
                R^{(2;k,j)}_{u;s,t}
            }{
                \omega(s,t)^{
                    \left(
                        \pfloor - j
                    \right)
                    /p
                }
            }
        \right|
    < \infty
    .
\end{align*}

The following series of results are easily verifiable by taking the same arguments as in the case where $X = \Xtilde$ and making suitable substitutions where needed.

\begin{lemma}\label{remainders are controlled paths lemma 2}
Let $R^{(2;k,j)}: [0,T]^2 \times [0, T] \rightarrow \Hom(V^{\otimes j}, \Hom(\Vtilde^{\otimes k}, E))$ and $R^{(1;j,k)}: [0,T]^2 \times [0, T] \rightarrow \Hom(\Vtilde^{\otimes k}, \Hom(V^{\otimes j}, E))$ be the maps defined by
\begin{align*}
    R^{(2;k,j)}_{u,v;s}(x)(y)
    &:=
    R^{(1;j,k)}_{s;u,v}(y)(x),
    \\
    R^{(1;j,k)}_{s,t;u}(x)(y)
    &:=
    R^{(2;k,j)}_{u;s,t}(y)(x).
\end{align*}
Then for any $0 \leq u \leq v \leq T$ and any $k = 0, \dotsc, \Ntilde$, the remainder tuple $\left(R^{(2;k,0)}_{u,v}, \dotsc, R^{(2;k,N)}_{u,v}\right)$ is an $X$-controlled path with remainders $\{\Rbf^{(2; k, j)}_{u,v} : j = 0, \dotsc, N \}$. Similarly for any $0 \leq s \leq v \leq T$, the tuples $\left( R^{(1;j,0)}_{s,t}, \dotsc, R^{(1;j,\Ntilde)}_{s,t}\right)$ for $j = 0, \dotsc, N$ are $\Xtilde$-controlled paths with remainders $\{\Rbf^{(1; j, k)}_{s,t} : k = 0, \dotsc, \Ntilde \}$.

These second order remainders are such that for any $x \in V^{\otimes j}$ and any $y \in \Vtilde^{\otimes k}$,
\begin{equation}
    \label{remainder relation 2}
    \Rbf^{(1; j, k)}_{s,t;u,v}
    (y)(x)
    =
    \Rbf^{(2; k, j)}_{u,v;s,t}
    (x)(y)
    .
\end{equation}
\end{lemma}

In place of Condition \ref{mixed variation condition} we now substitute with the following condition:

\begin{condition}\label{mixed variation condition 2}
Suppose that $\{Y^{(i;j,k)}\} \in \Dscr^{p, \ptilde}_{X, \Xtilde}([0,T]^2; E)$ with remainders $R^{(i;j,k)}$. For $j = 0, \dotsc, N$ and $k = 0, \dotsc, \Ntilde$, let $p_j = p/(j+1)$ and $\ptilde_k = \ptilde/(k+1)$. Assume there exists some $q_j$, $\qtilde_k$, such that
\begin{align*}
    \frac{1}{p_j}
    +
    \frac{1}{q_j}
    = \theta_j > 1
    ,\qquad
    \frac{1}{\ptilde_k}
    +
    \frac{1}{\qtilde_k}
    =
    \thetatil_k > 1
    ,
\end{align*}
and remainders $R^{(1;j,k)}_s$ have finite $\qtilde_k$-variation, $R^{(2;k,j)}_u$ have finite $q_j$-variation, and $\Rbf^{(1;j,k)}$ have finite $(\omega, \omegatil)$-controlled $(q_j, \qtilde_k)$-variation, by which we mean
\begin{align*}
    \left\|
        \Rbf^{(1;j,k)}
    \right\|_{(q_j, \omega), (\qtilde_k, \omegatil)}
    :=
    \sup_{\substack{
            0 \leq s < t \leq T
            \\
            0 \leq u < v \leq T
        }
    }
    \left|
        \frac{
            \Rbf^{(1;j,k)}_{s,t;u,v}
        }{
            \omega(s,t)^{1/q_j}
            \omegatil(u,v)^{1/\qtilde_k}
        }
    \right|
    .
\end{align*}
Define $\theta_* = \min_{j,k} \{\theta_j, \thetatil_k\}$ and $\theta^* = \max_{j,k} \{\theta_j, \thetatil_k\}$.
\end{condition}

With this condition in place we are again able to show that the one parameter rough integrals are controlled paths and thus can be integrated again.

\begin{lemma}
Suppose that Condition \ref{mixed variation condition 2} holds. For $j = 0, \dotsc, N$, $k = 0, \dotsc, \Ntilde$, define $Z^{(j)}_r$ and $\Ztilde^{(k)}_{r'}$ to be the following rough integrals for $r \in [s,t] \subset [0,T]$ and $r' \in [u,v] \subset [0,T]$,
\[
    Z_{r}^{(1;j)}
    =
    \int^v_u
        Y_{r, r'}^{(1;j,0)}
    d\Xtilde_{r'}
    ,\qquad
    Z_{r'}^{(2;k)}
    =
    \int^t_s
        Y_{r, r'}^{(2;k,0)}
    dX_{r}.
\]
The tuple $\left( Z_r^{(1;0)}, \dotsc, Z_r^{(1;N)} \right)$ is an $X$-controlled path and $\left( Z_{r'}^{(2;0)}, \dotsc, Z_{r'}^{(2;\Ntilde)} \right)$ is an $\Xtilde$-controlled path.
\end{lemma}

We again define local approximations of our joint integral and are able to show that these satisfy some identities involving remainders of $Y$.

\begin{lemma}
Define the local approximations
\[
    \Omega^Y
    \binom{s,t}{u,v}
    :=
    \sum_{j=0}^{N}
    \sum_{k=0}^{\Ntilde}
        Y_{s,u}^{(1;j,k)}
        \Big(
            \Xtilde^{k+1}_{u,v}
        \Big)
        \Big(
            X^{j+1}_{s,t}
        \Big)
    =
    \sum_{j=0}^{N}
    \sum_{k=0}^{\Ntilde}
        Y_{s,u}^{(2;k,j)}
        \Big(
            X^{j+1}_{s,t}
        \Big)
        \Big(
            \Xtilde^{k+1}_{u,v}
        \Big)
\]
with related quantities $\Gamma^Y$ and $\Theta^Y$ defined in the same way as in the $\Xtilde = X$ case (see Definition \ref{local approximation definitions}). These quantities satisfy the identities
\begin{align*}
    \Gamma^Y
    \binom{s,s',t}{u,v}
    &=
    \sum_{j=0}^{N}
    \sum_{k=0}^{\Ntilde}
        R^{(2;k,j)}_{u;s,s'}
        \Big(
            X^{j+1}_{s',t}
        \Big)
        \Big(
            \Xtilde^{k+1}_{u,v}
        \Big)
    \\[0.5em]
    \Theta^Y
    \binom{s,s',t}{u,u',v}
    &=
    \sum_{j=0}^{N}
    \sum_{k=0}^{\Ntilde}
        \Rbf^{(1;j,k)}_{s,s';u,u'}
        \Big(
            \Xtilde^{k+1}_{u',v}
        \Big)
        \Big(
            X^{j+1}_{s',t}
        \Big).
\end{align*}
\end{lemma}

Again, we repeat the same strategy as before to establish a maximal inequality, whereby we build a series of intermediate bounds using the above identity. Using the same type of calculations as in the $X = \Xtilde$ case, we are lead to the same type of results with $\Vbf^{(i;p)}_{[s,t], [u,v]}(X,Y)$ and $\eta^{(i;p,\alpha)}_{[s,t],[u,v]} (X,Y)$ replaced with the quantities
\begin{align*}
    \Vbf^{(1;p, \ptilde)}_{[s,t], [u,v]}(X, \Xtilde, Y)
    &:=
    \max_{
        \substack{
            j = 0, \dotsc, N
            \\
            k = 0, \dotsc, \Ntilde
        }
    }
        \left\{
            V^{p/(j+1)}_{[s,t]}\left(
                X^{j+1}
            \right)
            V^{\ptilde/(k+1)}_{[u,v]}\left(
                \Xtilde^{k+1}
            \right)
            V^{q_j, \qtilde_k}_{[s,t]\times[u,v]}\left(
                \Rbf^{(1; j, k)}
            \right)
        \right\}
    \\[0.5em]
    \Vbf^{(2;p, \ptilde)}_{[s,t], [u,v]}(X, \Xtilde, Y)
    &:=
    \max_{
        \substack{
            j = 0, \dotsc, N
            \\
            k = 0, \dotsc, \Ntilde
        }
    }
        \left\{
            V^{p/(j+1)}_{[s,t]}\left(
                X^{j+1}
            \right)
            V^{\ptilde/(k+1)}_{[u,v]}\left(
                \Xtilde^{k+1}
            \right)
            V^{\qtilde_k, q_j}_{[s,t]\times[u,v]}\left(
                \Rbf^{(2; k, j)}
            \right)
        \right\}
    \\[0.5em]
    \eta^{(1;p,\ptilde,\alpha)}_{[s,t],[u,v]}
    (X, \Xtilde, Y)
    &:=
    \left(
        \sum^{N}_{j= 0}
        \sum^{\Ntilde}_{k = 0}
            \left|
                \Xtilde^{k+1}_{u, v}
            \right|^{\alpha}
            V^{p/(j+1)}_{[s,t]}
            \left(
                X^{j+1}
            \right)^{\alpha}
            V^{q_j}_{[s,t]}
            \left(
                R^{(2;k,j)}_u
            \right)^{\alpha}
    \right)^{1/\alpha}
    \\[0.5em]
    \eta^{(2;p,\ptilde,\alpha)}_{[s,t],[u,v]}
    (X, \Xtilde, Y)
    &:=
    \left(
        \sum^{N}_{j= 0}
        \sum^{\Ntilde}_{k = 0}
            \left|
                X^{j+1}_{s, t}
            \right|^{\alpha}
            V^{\ptilde/(k+1)}_{[u,v]}
            \left(
                \Xtilde^{k+1}
            \right)^{\alpha}
            V^{\qtilde_k}_{[u,v]}
            \left(
                R^{(1;j,k)}_s
            \right)^{\alpha}
    \right)^{1/\alpha}.
\end{align*}

Making these substitutions into the intermediate bounds of Lemmas \ref{Theta inequality lemma}, \ref{Delta inequalities lemma}, and \ref{endpoints lemma}, the maximal inequality then takes the following form:

\begin{theorem}[Maximal inequality]
Let Condition \ref{mixed variation condition 2} hold. There exists a constant $C(p, \ptilde \alpha, \theta_*)$ such that for all partitions $\Dcal \subset [s,t]$ and $\Dcal' \subset [u,v]$ the following inequality holds
\begin{equation*}
    \left|
        \sum_{\Dcal \times \Dcal'}
            \Omega^Y
        -
        \Omega^Y
        \binom{s,t}{u,v}
    \right|
    \leq
    C(p, \ptilde, \alpha, \theta_*)
    \left(
        \Vbf^{(p, \ptilde)}_{[s,t], [u,v]}
        (X, \Xtilde, Y)
        +
        \eta^{(p, \ptilde, \alpha)}_{[s,t],[u,v]}
        (X, \Xtilde, Y)
    \right)
\end{equation*}
where
\begin{align*}
    \Vbf^{(p, \ptilde)}_{[s,t], [u,v]}
    (X, \Xtilde, Y)
    &=
    \min\left(
        \Vbf^{(1;p, \ptilde)}_{[s,t], [u,v]}(X, \Xtilde, Y)
        ,
        \Vbf^{(2;p, \ptilde)}_{[s,t], [u,v]}(X, \Xtilde, Y)
    \right)
    \\[0.5em]
    \eta^{(p, \ptilde, \alpha)}_{[s,t],[u,v]}
    (X, \Xtilde, Y)
    &=
    \eta^{(1;p,\ptilde,\alpha)}_{[s,t],[u,v]}
    (X, \Xtilde, Y)
    +
    \eta^{(2;p,\ptilde,\alpha)}_{[s,t],[u,v]}
    (X, \Xtilde, Y).
\end{align*}
\end{theorem}

Once the maximal inequality is established, existence of the joint integral, a rough Fubini type theorem, and a stability result can be formulated, under some additional uniformity and regularity conditions.

\begin{theorem}[Existence of integral and a rough Fubini type theorem]
\label{integral existence theorem 2}
Assume that Condition \ref{mixed variation condition 2} holds and also assume the additional uniformity condition on remainders
\[
    \sup_{0 \leq s \leq T}
        \left\|
            R^{(1;j,k)}_{s}
        \right\|_{\qtilde_k, \omegatil}
    \leq
    C
    , \qquad
    \sup_{0 \leq u \leq T}
        \left\|
            R^{(2;k,j)}_{u}
        \right\|_{q_j, \omega}
    \leq
    C
    ,
\]
for some constant $C$.

Then we have existence of the integral
\begin{equation}\label{integral existence 2}
    \int_{[s,t] \times [u,v]}
        Y_{r,r'}
    \;d(X_r, \Xtilde_{r'})
    :=
    \lim_{|\Dcal \times \Dcal'| \rightarrow 0}
        \sum_{\Dcal \times \Dcal'}
            \Omega^Y
\end{equation}
and we have the bound
\begin{equation}\label{integral bound 2}
    \left|
        \int_{[s,t] \times [u,v]}
            Y_{r,r'}
        \;d(X_r, \Xtilde_{r'})
        -
        \Omega^Y
        \binom{s,t}{u,v}
    \right|
    \leq
    C_{p, \ptilde}
    \left(
        \Vbf^{(p, \ptilde)}_{[s,t], [u,v]}
        (X, \Xtilde, Y)
        +
        \eta^{(p, \ptilde, \alpha)}_{[s,t],[u,v]}
        (X, \Xtilde, Y)
    \right),
\end{equation}
for some constant $C_{p, \ptilde} = C(p, \ptilde, \alpha, \theta_*)$. Moreover the joint integral is equal to the double integrals,
\begin{equation}\label{rough Fubini 2}
    \int_s^t
        \left(
            \int_u^v
                Y_{r,r'}
            \;d\Xtilde_{r'}
        \right)
    \;dX_{r}
    =
    \int_{[s,t] \times [u,v]}
        Y_{r,r'}
    \;d(X_r, \Xtilde_{r'})
    =
    \int_u^v
        \left(
            \int_s^t
                Y_{r,r'}
            \;dX_{r}
        \right)
    \;d\Xtilde_{r'}
    .
\end{equation}
\end{theorem}

\begin{theorem}[Stability of double integral]
For $i = 1,2$, let $X^{(i)} = (1, X^{(i;1)}, \dotsc, X^{(i;\pfloor)})$ be two geometric $p$-rough paths on $V$ with control $\omega$ and similarly consider two geometric $\ptilde$-rough paths $\Xtilde^{(i)} = (i, \Xtilde^{(i;1)}, \dotsc, \Xtilde^{(i;\ptilfloor)})$ on $\Vtilde$ with control $\omegatil$. Denote by $\Xbf^{(i)}$ the pairs $\Xbf^{(i)} = (X^{(i)}, \Xtilde^{(i)})$.

Consider two jointly controlled paths $\{Y^{(i;j,k)}\} \in \Dscr_{\Xbf^{(1)}}^{p, \ptilde}$ and $\{\Ytilde^{(i;j,k)}\} \in \Dscr_{\Xbf^{(2)}}^{p, \ptilde}$ with first order remainders $\{R^{(i;j,k)}\}, \{\Rtilde^{(i;j,k)}\}$ and second order remainders $\{\Rbf^{(i;j,k)}\}, \{\Rbftilde^{(i;j,k)}\}$ respectively. Fix $[s,t], [u,v] \subset [0,T]$ and suppose that the two paths satisfy the conditions of Theorem \ref{integral existence theorem 2}, and define the following "distance" between two jointly controlled paths
\begin{alignat*}{2}
    d_{\Xbf^{(1)}, \Xbf^{(2)}}
    (Y, \Ytilde)
    &:=
    \sum^{N}_{j=0}
    \sum^{\Ntilde}_{k=0}
    &&
        \left|
            Y^{(1;j,k)}_{s,u}
            -
            \Ytilde^{(1;j,k)}_{s,u}
        \right|
        +
        \left\|
            R^{(1;j,k)}_{s}
            -
            \Rtilde^{(1;j,k)}_{s}
        \right\|_{\qtilde_k, \omegatil}
    \\
    &
    &&
        +
        \left\|
            R^{(2;k,j)}_{u}
            -
            \Rtilde^{(2;k,j)}_{u}
        \right\|_{q_j, \omega}
        +
        \left\|
            \Rbf^{(1;j,k)}
            -
            \Rbftilde^{(1;j,k)}
        \right\|_{(q_j, \omega), (\qtilde_k, \omegatil)}
    .
\end{alignat*}
The double integrals
\[
    I_1
    =
    \int_{[s,t] \times [u,v]}
        Y_{r,r'}
    \;d
    \left(
        X^{(1)}_r, \Xtilde^{(1)}_{r'}
    \right)
    ,\qquad
    I_2
    =
    \int_{[s,t] \times [u,v]}
        \Ytilde_{r,r'}
    \;d
    \left(
        X_r^{(2)}, \Xtilde_{r'}^{(2)}
    \right)
    ,
\]
are such that
\begin{equation*}
    \left|
        I_1
        -
        I_2
    \right|
    \leq
    C
    \left(
        \|X^{(1)} - X^{(2)}\|_{p\variation}
        +
        \|\Xtilde^{(1)} - \Xtilde^{(2)}\|_{\ptilde\variation}
        +
        d_{\Xbf^{(1)}, \Xbf^{(2)}}(Y, \Ytilde)
    \right)
    ,
\end{equation*}
where $C$ depends on $\Xbf^{(1)}, \Xbf^{(2)}, Y, \Ytilde$, and remainders, at initial and terminal times $s,t,u,v$.
\end{theorem}

To finish we give an example of a jointly controlled path satisfying the conditions of Theorem \ref{integral existence theorem 2} in the form of the signature kernel, which has recently seen applications in data science in \cite{sigkernelsgoursat}. The signature kernel $K_{(s,t), (u,v)}(X, \Xtilde) = \langle S(X)_{s,t}, S(\Xtilde)_{u,v} \rangle$ also happens to satisfy a two-parameter rough integral equation as can be seen in Theorem 4.11 in \cite{sigkernelsgoursat}.

\begin{lemma}
Consider two $\omega$-controlled geometric $p$-rough paths $X, \Xtilde$ on $V$. Fix $u_0, s_0 \in [0, T]$. Define the path $Y:[s_0, T] \times [u_0, T] \rightarrow \Real$ by
\begin{equation*}
    Y_{s,u}
    =
    K_{(s_0, s), (u_0, u)}
    (X, \Xtilde)
    =
    \left\langle
        S(X)_{s_0, s}
        ,
        S(\Xtilde)_{u_0, u}
    \right\rangle_{\overline{T(V)}}
    =
    \sum_{l=0}^{\infty}
        \left\langle
            X^{l}_{s_0, s}
            ,
            \Xtilde^{l}_{u_0, u}
        \right\rangle_{V^{\otimes l}}
\end{equation*}
The path $Y$ is a jointly $(X, \Xtilde)$-controlled path with Gubinelli derivatives
\begin{align}
    \label{Gubinelli derivatives of the signature kernel}
    Y^{(1;j,k)}_{s,u}
    (y)(x)
    =
    Y^{(2;k,j)}_{s,u}
    (x)(y)
    =
    \begin{cases}
        \displaystyle\sum^{\infty}_{l=j}
            \left\langle
                \Xtilde_{u_0, u}^{l-k}
                \otimes
                y
                ,\;
                X^{l-j}_{s_0, s}
                \otimes
                x
            \right\rangle_{V^{\otimes l}}
            ,
        &
        k \leq j,
        \\[1.2em]
        \displaystyle\sum^{\infty}_{l=k}
            \left\langle
                \Xtilde_{u_0, u}^{l-k}
                \otimes
                y
                ,\;
                X^{l-j}_{s_0, s}
                \otimes
                x
            \right\rangle_{V^{\otimes l}}
            ,
        &
        k \geq j,
    \end{cases}
\end{align}
and remainders
\begin{align}
    \label{signature kernel remainder 1}
    R^{(1;j,k)}_{s;u,v}
    (y)(x)
    &=
    \sum_{n = \pfloor - k}^{\infty}
    \;
    \sum_{l = m + k}^{\infty}
        \left\langle
            \Xtilde^{l-k-n}_{u_0, u}
            \otimes
            \Xtilde^{n}_{u,v}
            \otimes
            y
            ,\;
            X^{l-j}_{s_0, s}
            \otimes
            x
        \right\rangle_{V^{\otimes l}}
    ,
    \\[0.5em]
    \label{signature kernel remainder 2}
    R^{(2;k,j)}_{u;s,t}
    (x)(y)
    &=
    \sum_{m = \pfloor - j}^{\infty}
    \;
    \sum_{l = m + j}^{\infty}
        \left\langle
            X^{l-j-m}_{s_0, s}
            \otimes
            X^{m}_{s,t}
            \otimes
            x
            ,\;
            \Xtilde^{l-k}_{u_0, u}
            \otimes
            y
        \right\rangle_{V^{\otimes l}}
    ,
    \\[0.5em]
    \label{second order kernel remainder}
    \Rbf^{(1;j,k)}_{s,t;u,v}
    (y)(x)
    &=
    \sum_{l = \pfloor}^{\infty}
    \;
    \sum_{m = \pfloor - j}^{l-j}
    \;
    \sum_{n = \pfloor - k}^{l-k}
        \left\langle
            X^{l-j-m}_{s_0, s}
            \otimes
            X^{m}_{s, t}
            \otimes
            x
            ,\;
            \Xtilde^{l-k-n}_{u_0, u}
            \otimes
            \Xtilde^{n}_{u, v}
            \otimes
            y
        \right\rangle_{V^{\otimes l}}
    .
\end{align}
\end{lemma}

\begin{proof}
Because of the symmetry of the signature kernel, showing that $\left( Y^{(2;k,0)}_{\cdot,u}, \dotsc, Y^{(2;k,N)}_{\cdot,u} \right)$ is an $X$-controlled path follows very similarly to showing $\left( Y^{(1;j,0)}_{s,\cdot}, \dotsc, Y^{(1;j,N)}_{u, \cdot} \right)$ is an $\Xtilde$-controlled path, so we only show the latter. Upon inspection the remainders satisfy the required regularity conditions since
\[
    \left\|
        X^l_{s,t}
    \right\|_{V^{\otimes l}}
    \leq
    \frac{
        \omega(s,t)^{l/p}
    }{
        \beta_p (l/p)!
    }
    ,
    \qquad
    \left\|
        \Xtilde^l_{u,v}
    \right\|_{V^{\otimes l}}
    \leq
    \frac{
        \omega(u,v)^{l/p}
    }{
        \beta_p (l/p)!
    }
\]
where $(l/p)! = \Gamma(l/p)$ and $\beta_p$ is a constant which depends only on $p$. This also gives us that the second order remainders have the appropriate finite $\omega$-controlled mixed variation required in Theorem \ref{integral existence theorem 2}.

We begin by looking for suitable candidates for $Y^{(2;0,j)}$ in order to determine $Y^{(1;j,0)}$. Rewriting the signature kernel with Chen's identity, we have
\begin{align*}
    Y_{t,u}
    &=
    \sum^{\infty}_{l=0}
        \left\langle
            \sum^{l}_{j=0}
                X^{l-j}_{s_0, s}
                \otimes
                X^{j}_{s,t}
            ,\;
            \Xtilde^{l}_{u_0, u}
        \right\rangle_{V^{\otimes l}}
    \\[0.5em]
    &=
    \sum^{N}_{j=0}
    \sum^{\infty}_{l=j}
        \left\langle
            X^{l-j}_{s_0, s}
            \otimes
            X^{j}_{s,t}
            ,\;
            \Xtilde^{l}_{u_0, u}
        \right\rangle_{V^{\otimes l}}
    +
    \sum^{\infty}_{j=\pfloor}
    \sum^{\infty}_{l=j}
        \left\langle
            X^{l-j}_{s_0, s}
            \otimes
            X^{j}_{s,t}
            ,\;
            \Xtilde^{l}_{u_0, u}
        \right\rangle_{V^{\otimes l}}
    .
\end{align*}
From this we deduce that we want to have
\[
    Y^{(2;0,j)}_{s,u}
    (x)
    =
    \sum^{\infty}_{l=j}
        \left\langle
            X^{l-j}_{s_0, s}
            \otimes
            x
            ,\;
            \Xtilde^{l}_{u_0, u}
        \right\rangle_{V^{\otimes l}}
    ,
\]
and further calculations yield
\begin{align*}
    Y^{(2;0,j)}_{t,u}
    (x)
    &=
    \sum^{\infty}_{l=j}
        \left\langle
            \sum_{m=0}^{l-j}
                X^{l-j-m}_{s_0, s}
                \otimes
                X^{m}_{s,t}
                \otimes
                x
            ,\;
            \Xtilde^{l}_{u_0, u}
        \right\rangle_{V^{\otimes l}}
    \\[0.5em]
    &=
    \sum_{m=0}^{\infty}
    \;
    \sum^{\infty}_{l=j+m}
        \left\langle
            X^{l-j-m}_{s_0, s}
            \otimes
            X^{m}_{s,t}
            \otimes
            x
            ,\;
            \Xtilde^{l}_{u_0, u}
        \right\rangle_{V^{\otimes l}}
    \\[0.5em]
    &=
    \sum_{m=0}^{N-j}
        Y^{(2;0,j)}_{s,u}
        \left(
            X^{m}_{s,t} \otimes x
        \right)
    +
    \sum_{m=\pfloor -j}^{\infty}
    \;
    \sum^{\infty}_{l=j+m}
        \left\langle
            X^{l-j-m}_{s_0, s}
            \otimes
            X^{m}_{s,t}
            \otimes
            x
            ,\;
            \Xtilde^{l}_{u_0, u}
        \right\rangle_{V^{\otimes l}}
    .
\end{align*}
We now perform a similar procedure on $Y^{(1;j,0)}$ to obtain the following:
\begin{align*}
    Y^{(1;j,0)}_{s,v}
    (x)
    &=
    \sum^{\infty}_{l=j}
        \left\langle
            \sum^{l}_{n=0}
                \Xtilde^{l-n}_{u_0, u}
                \otimes
                \Xtilde^{n}_{u, v}
            ,\;
            X^{l-j}_{s_0, s}
            \otimes
            x
        \right\rangle_{V^{\otimes l}}
    \\[0.5em]
    &=
    \sum_{n=0}^{j}
    \;
    \sum^{\infty}_{l=j}
    \left\langle
            \Xtilde^{l-n}_{u_0, u}
            \otimes
            \Xtilde^{n}_{u, v}
            ,\;
            X^{l-j}_{s_0, s}
            \otimes
            x
        \right\rangle_{V^{\otimes l}}
    +
    \sum_{n=j+1}^{\infty}
    \;
    \sum^{\infty}_{l=j}
    \left\langle
            \Xtilde^{l-n}_{u_0, u}
            \otimes
            \Xtilde^{n}_{u, v}
            ,\;
            X^{l-j}_{s_0, s}
            \otimes
            x
        \right\rangle_{V^{\otimes l}}
    ,
\end{align*}
which leads us to choose $Y^{(1;j,k)}$ as in \eqref{Gubinelli derivatives of the signature kernel}. The final step is now to verify that these derivatives satisfy the properties of controlled paths. We split into two cases, first where $k \leq j$,
\begin{align*}
    Y^{(1;j,k)}_{s,v}
    (y)(x)
    &=
    \sum^{\infty}_{l=j}
        \left\langle
            \sum_{n=0}^{l-k}
                \Xtilde_{u_0, u}^{l-k-n}
                \otimes
                \Xtilde_{u, v}^{n}
                \otimes
                y
            ,\;
            X^{l-j}_{s_0, s}
            \otimes
            x
        \right\rangle_{V^{\otimes l}}
    \\[0.5em]
    &=
    \sum_{n=0}^{j-k}
    \;
    \sum^{\infty}_{l=j}
        \left\langle
            \Xtilde_{u_0, u}^{l-k-n}
            \otimes
            \Xtilde_{u, v}^{n}
            \otimes
            y
            ,\;
            X^{l-j}_{s_0, s}
            \otimes
            x
        \right\rangle_{V^{\otimes l}}
    \\
    &\quad
    +
    \sum_{n=j-k+1}^{\infty}
    \;
    \sum^{\infty}_{l=n+k}
        \left\langle
            \Xtilde_{u_0, u}^{l-k-n}
            \otimes
            \Xtilde_{u, v}^{n}
            \otimes
            y
            ,\;
            X^{l-j}_{s_0, s}
            \otimes
            x
        \right\rangle_{V^{\otimes l}}
    \\[0.5em]
    &=
    \sum_{n=0}^{N-k}
        Y^{(1;j,k+n)}_{s,u}
        \left(
            \Xtilde^n_{u,v}
            \otimes
            y
        \right)
        (x)
    \\
    &\quad
    +
    \sum_{n = \pfloor -k}^{\infty}
    \;
    \sum^{\infty}_{l=n+k}
        \left\langle
            \Xtilde_{u_0, u}^{l-k-n}
            \otimes
            \Xtilde_{u, v}^{n}
            \otimes
            y
            ,\;
            X^{l-j}_{s_0, s}
            \otimes
            x
        \right\rangle_{V^{\otimes l}}
    ,
\end{align*}
and in the case $j \leq k$ we have
\begin{align*}
    Y^{(1;j,k)}_{s,v}
    (y)(x)
    &=
    \sum^{\infty}_{l=k}
        \left\langle
            \sum_{n=0}^{l-k}
                \Xtilde_{u_0, u}^{l-k-n}
                \otimes
                \Xtilde_{u, v}^{n}
                \otimes
                y
            ,\;
            X^{l-j}_{s_0, s}
            \otimes
            x
        \right\rangle_{V^{\otimes l}}
    \\[0.5em]
    &=
    \sum_{n = 0}^{\infty}
    \;
    \sum^{\infty}_{l=n+k}
        \left\langle
            \Xtilde_{u_0, u}^{l-k-n}
            \otimes
            \Xtilde_{u, v}^{n}
            \otimes
            y
            ,\;
            X^{l-j}_{s_0, s}
            \otimes
            x
        \right\rangle_{V^{\otimes l}}
    \\[0.5em]
    &=
    \sum_{n=0}^{N-k}
        Y^{(1;j,k+n)}_{s,u}
        \left(
            \Xtilde^n_{u,v}
            \otimes
            y
        \right)
        (x)
    \\
    &\quad
    +
    \sum_{n = \pfloor -k}^{\infty}
    \;
    \sum^{\infty}_{l=n+k}
        \left\langle
            \Xtilde_{u_0, u}^{l-k-n}
            \otimes
            \Xtilde_{u, v}^{n}
            \otimes
            y
            ,\;
            X^{l-j}_{s_0, s}
            \otimes
            x
        \right\rangle_{V^{\otimes l}}
    .
\end{align*}
By definition, for the second order remainder we have
\begin{align*}
    \Rbf^{(1;j,k)}_{s,t;u,v}
    (y)(x)
    &=
    R^{(1;j,k)}_{t;u,v}
    -
    \sum_{m=0}^{N-j}
        R^{(1;j+m,k)}_{s;u,v}
        (y)
        \left(
           X^{m}_{s,t}
        \right)
        (x)
    \\[0.5em]
    &=
    \sum_{n = \pfloor -k}^{\infty}
    \;
    \sum^{\infty}_{l=n+k}
        \left\langle
            \Xtilde_{u_0, u}^{l-k-n}
            \otimes
            \Xtilde_{u, v}^{n}
            \otimes
            y
            ,\;
            X^{l-j}_{s_0, t}
            \otimes
            x
            -
            \sum_{m=0}^{N-j}
                X^{l-j-m}_{s_0, s}
                \otimes
                X^{m}_{s,t}
                \otimes
                x
        \right\rangle_{V^{\otimes l}}
    \\[0.5em]
    &=
    \sum_{n = \pfloor -k}^{\infty}
    \;
    \sum^{\infty}_{l=n+k}
        \left\langle
            \Xtilde_{u_0, u}^{l-k-n}
            \otimes
            \Xtilde_{u, v}^{n}
            \otimes
            y
            ,\;
            \sum_{m=\pfloor - j}^{l-j}
                X^{l-j-m}_{s_0, s}
                \otimes
                X^{m}_{s,t}
                \otimes
                x
        \right\rangle_{V^{\otimes l}}
    \\[0.5em]
    &=
    \sum^{\infty}_{l = \pfloor}
    \;
    \sum_{m = \pfloor - j}^{l-j}
    \;
    \sum_{n = \pfloor -k}^{l-k}
        \left\langle
            \Xtilde_{u_0, u}^{l-k-n}
            \otimes
            \Xtilde_{u, v}^{n}
            \otimes
            y
            ,\;
            X^{l-j-m}_{s_0, s}
            \otimes
            X^{m}_{s,t}
            \otimes
            x
        \right\rangle_{V^{\otimes l}}
    .
\end{align*}
\end{proof}


\printbibliography

@article{Lyo98,
    author = {Lyons, Terry J.},
    journal = {Revista Matemática Iberoamericana},
    number = {2},
    pages = {215-310},
    title = {Differential equations driven by rough signals.},
    url = {http://eudml.org/doc/39555},
    volume = {14},
    year = {1998},
}

@article{Gub04,
    title = {Controlling rough paths},
    journal = {Journal of Functional Analysis},
    volume = {216},
    number = {1},
    pages = {86-140},
    year = {2004},
    issn = {0022-1236},
    doi = {https://doi.org/10.1016/j.jfa.2004.01.002},
    url = {https://www.sciencedirect.com/science/article/pii/S0022123604000497},
    author = {Gubinelli, Massimiliano},
}

@article{Gub10,
    title = {Ramification of rough paths},
    journal = {Journal of Differential Equations},
    volume = {248},
    number = {4},
    pages = {693-721},
    year = {2010},
    issn = {0022-0396},
    doi = {https://doi.org/10.1016/j.jde.2009.11.015},
    url = {https://www.sciencedirect.com/science/article/pii/S0022039609004379},
    author = {Gubinelli, Massimiliano},
}

@article{FdLP06,
    author = {Denis Feyel and Arnaud de La Pradelle},
    title = {{Curvilinear Integrals Along Enriched Paths}},
    volume = {11},
    journal = {Electronic Journal of Probability},
    number = {none},
    publisher = {Institute of Mathematical Statistics and Bernoulli Society},
    pages = {860 -- 892},
    year = {2006},
    doi = {10.1214/EJP.v11-356},
    URL = {https://doi.org/10.1214/EJP.v11-356}
}

@book{CLL07,
    author = {Caruana, Michael J. and L\'evy, Thierry and Lyons, Terry, J.},
    year = {2007},
    title = {Differential Equations Driven by Rough Paths},
    series = {Lecture Notes in Mathematics},
    edition = {1},
    publisher = {Springer, Berlin, Heidelberg},
    doi = {10.1007/978-3-540-71285-5},
}

@book{FH14,
    author = {Friz, Peter K. and Hairer, Martin},
    year = {2014},
    title = {{A Course on Rough Paths: With an Introduction to Regularity Structures}},
    edition = {1},
    publisher = {Springer International Publishing},
    doi = {10.1007/978-3-319-08332-2},
}

@article{Tow02,
    author = {Towghi, Nasser},
    journal = {JIPAM. Journal of Inequalities in Pure \& Applied Mathematics [electronic only]},
    number = {2},
    pages = {Paper No. 22, 13 p., electronic only-Paper No. 22, 13 p., electronic only},
    publisher = {Victoria University, School of Communications and Informatics},
    title = {Multidimensional extension of L. C. Young's inequality.},
    url = {http://eudml.org/doc/122136},
    volume = {3},
    year = {2002},
}

@article{FV11,
    author = {Friz, Peter K. and Victoir, Nicholas},
    title = {{A Note on Higher Dimensional p-Variation}},
    volume = {16},
    journal = {Electronic Journal of Probability},
    number = {none},
    publisher = {Institute of Mathematical Statistics and Bernoulli Society},
    pages = {1880 -- 1899},
    year = {2011},
    doi = {10.1214/EJP.v16-951},
    URL = {https://doi.org/10.1214/EJP.v16-951}
}

@article{GH19,
    author = {Gerasimovičs, Andris and Hairer, Martin},
    title = {{Hörmander’s theorem for semilinear SPDEs}},
    volume = {24},
    journal = {Electronic Journal of Probability},
    number = {none},
    publisher = {Institute of Mathematical Statistics and Bernoulli Society},
    pages = {1 -- 56},
    year = {2019},
    doi = {10.1214/19-EJP387},
    URL = {https://doi.org/10.1214/19-EJP387}
}

@misc{CG14,
      title = {Rough sheets}, 
      author = {K. Chouk and M. Gubinelli},
      year = {2014},
      eprint = {1406.7748},
      archivePrefix = {arXiv},
      primaryClass = {math.PR}
}

@misc{sigkernelsgoursat,
    title={The Signature Kernel is the solution of a Goursat PDE}, 
    author={Salvi, Cristopher and Cass, Thomas and Foster, James and Lyons, Terry and Yang, Weixin},
    year={2021},
    eprint={2006.14794},
    archivePrefix={arXiv},
    primaryClass={math.AP},
    note={Accepted to a peer in SIAM Journal on Mathematics of Data Science},
}

@article{You36,
    author = {L. C. Young},
    title = {{An inequality of the Hölder type, connected with Stieltjes integration}},
    volume = {67},
    journal = {Acta Mathematica},
    number = {none},
    publisher = {Institut Mittag-Leffler},
    pages = {251 -- 282},
    year = {1936},
    doi = {10.1007/BF02401743},
    URL = {https://doi.org/10.1007/BF02401743}
}

@article{FGGR16,
    author = {Peter K. Friz and Benjamin Gess and Archil Gulisashvili and Sebastian Riedel},
    title = {{The Jain– Monrad criterion for rough paths and applications to random Fourier series and non-Markovian Hörmander theory}},
    volume = {44},
    journal = {The Annals of Probability},
    number = {1},
    publisher = {Institute of Mathematical Statistics},
    pages = {684 -- 738},
    year = {2016},
    doi = {10.1214/14-AOP986},
    URL = {https://doi.org/10.1214/14-AOP986}
}

\appendix

\section{Rough integration of controlled paths}

Here we look at the one time variable integration theory of controlled paths and write some results with explicit bounds which we use throughout. The objective here is to familiarise the reader with the choices of notation used in the paper and to give some insight into the strategy used for constructing the rough integral. Again, for the less experienced reader we suggest the works \cite{Lyo98, CLL07, FH14} as comprehensive introductions to the theory of rough paths.

Denote by $\Delta_T$ the 2-simplex $\Delta_T = \{ (s,t) \mid 0 \leq s \leq t \leq T\}$ and $\Delta_T^3$ the 3-simplex $\Delta_T^3 = \{ (s,s',t) \mid 0 \leq s \leq s' \leq t \leq T\}$. For a control $\omega$, Banach space $E$, $p > 0$, and $A: \Delta_T \rightarrow E$, define the norm
\[
    \| A \|_{p, \omega}
    =
    \sup_{0 \leq s < t \leq T}
        \left|
            \frac{A_{s.t}}{\omega(s,t)^{1/p}}
        \right|
    ,
\]
and define as the space $C_{\omega}^{p}([0,T];E)$ as the space of additive functions with finite $\| \cdot \|_{\alpha, \omega}$-norms. Similarly for $\beta > 0$ and $B: \Delta^3_T \rightarrow E$, we make a slight abuse of notation and define
\[
    \| B \|_{\beta, \omega}
    =
    \sup_{0 \leq s < s' < t \leq T}
        \left|
            \frac{B_{s,s',t}}{\omega(s,t)^{1/\beta}}
        \right|
    .
\]
We will drop the $\omega$ from the subscript when the choice of control is clear from the context.

The notion of controlled paths was first introduced by Gubinelli in \cite{Gub04} in the case $p < 3$, and was extended to a more general framework involving branched rough paths in \cite{Gub10}. Here we define controlled paths in the case where $V$ is finite-dimensional and $X$ is a geometric rough path on $V$. The definition parallels the notion of a $\Lip(\gamma)$ function which is used in the rough integration of one-forms, c.f. \cite{Lyo98, CLL07}.

\begin{definition}[Controlled paths in one variable]
\label{controlled path definition}
Let $X = (1, X^1, X^2, \dotsc, X^{\pfloor})$ be an $\omega$-controlled geometric $p$-rough path on a Banach space $V$. Suppose that $Y : [0, T] \rightarrow E$ for some Banach space $E$. Suppose for $j = 0, 1, \dotsc, N = \pfloor - 1$ there exists $Y^{(j)} : [0, T] \rightarrow \Hom(V^{\otimes j}, E)$, with $Y^{(0)} = Y$, such that $Y^{(j)} \in C_{\omega}^{1/p}([0,T];E)$ (where we associate $Y^{(j)}$ with $Y^{(j)}_{s,t} = Y^{(j)}_t - Y^{(j)}_s$) and 
\begin{align*}
    Y^{(j)}_t
    =
    \sum^{N - j}_{k = 0}
        Y^{(j+k)}_s
        \left(
            X^{k}_{s,t}
        \right)
    +
    R^{(j)}_{s,t},
\end{align*}
for remainders $R^{(j)}$ satisfying the regularity condition $\|R^{(j)}\|_{p/(\pfloor-j)} < \infty$.

Then we say that the tuple $\left( Y, Y^{(1)}, \dotsc, Y^{(N)}\right)$ is an $X$-controlled path over $[0,T]$ with values in $E$. Denote this space of $X$-controlled paths by $\Dscr_X^p ([0,T]; E)$. For $j = 1, \dotsc, N$ we call $Y^{(j)}$ Gubinelli derivatives of $Y$.
\end{definition}

The definition of controlled paths lends itself naturally to considering enhanced Riemann sums; for a controlled path $Y$ consider the local approximation
\begin{equation}\label{local approximation}
    \Xi^Y_{s,t}
    =
    \sum^N_{j=0}
        Y^{(j)}_s
        \left(
            X^{j+1}_{s,t}
        \right)
    ,
\end{equation}
and the enhanced Riemann sum
\[
    \sum_{\Dcal}
        \Xi^Y
    :=
    \sum_{m = 1}^{m_0}
        \Xi^Y_{s_{m-1}, s_{m}}
\]
over a partition $\Dcal = \{s_0 < \dotsc < s_{m_0} \} \subset [s,t]$. For a map $\Xi: [0, T]^2 \rightarrow E$, we will use $\delta \Xi: [0,T]^3 \rightarrow E$ to denote the map
\begin{equation}\label{delta Xi}
    \delta \Xi_{s,s',t}
    =
    \Xi_{s,t}
    -
    \Xi_{s,s'}
    -
    \Xi_{s',t}
    ,
\end{equation}
which is a quantity that appears naturally when considering sums over partitions, such as with enhanced Riemann sums. For a controlled path $Y$, the approximations $\Xi^Y$ are such that $\delta\Xi^Y$ can be expressed in terms of the remainders of $Y$. This follows as a consequence of Chen's identity as seen below.

\begin{lemma}\label{controlled path identity lemma}
Let $(Y^{(0)}, \dotsc, Y^{(N)})$ be an $X$-controlled path with remainders $\left\{R^{(j)}\right\}_{j=0}^N$. Define the two parameter process $\Xi^Y$ by \eqref{local approximation}.

Then for $0 \leq s \leq s' \leq t \leq T$,
\begin{equation} \label{controlled path identity}
    -\delta \Xi^Y_{s,s',t}
    =
    \sum_{j=0}^N
        \left(
            Y_s^{(j)}
            \left(
                X^{j+1}_{s,s'}
            \right)
            +
            Y_{s'}^{(j)}
            \left(
                X^{j+1}_{s',t}
            \right)
            -
            Y_s^{(j)}
            \left(
                X^{j+1}_{s,t}
            \right)
        \right)
    =
    \sum_{j=0}^N
        R^{(j)}_{s,s'}
        \left(
            X^{j+1}_{s',t}
        \right)
    .
\end{equation}
\end{lemma}

\begin{proof}
We first note that by Chen's identity
\begin{equation*}
    X_{s,t}^{j+1}
    -
    X_{s,s'}^{j+1}
    -
    X_{s',t}^{j+1}
    =
    \sum_{l=1}^j
        X_{s,s'}^{j+1-l}
        \otimes
        X_{s',t}^{l}
    =
    \sum_{l=0}^{j-1}
        X_{s,s'}^{j-l}
        \otimes
        X_{s',t}^{l+1}.
\end{equation*}
Using this identity we then write
\begin{align*}
    -\delta \Xi^Y_{s,s',t}
    &=
    \sum_{j=0}^{N}
        \left(
            \left(
                Y_{s'}^{(j)}
                -
                Y_s^{(j)}
            \right)
            \left(
                X^{j+1}_{s',t}
            \right)
            +
            Y_s^{(j)}
            \left(
                X^{j+1}_{s,s'}
                +
                X^{j+1}_{s',t}
                -
                X^{j+1}_{s,t}
            \right)
        \right)
    \\[0.5em]
    &=
    \sum_{j=0}^{N}
        \left(
        \sum_{m=1}^{N-j}
            Y^{(j+m)}_{s}
            \left(
                X^{m}_{s,s'}
            \right)
            +
            R^{(j)}_{s,s'}
        \right)
        \left(
            X^{j+1}_{s',t}
        \right)
    -
    \sum_{j=1}^N
        \sum_{l=0}^{j-1}
            Y_s^{(j)}
            \Big(
                X_{s,s'}^{j-l}
            \Big)
            \Big(
                X_{s',t}^{l+1}
            \Big)
    \\[0.5em]
    &=
    \sum_{j=0}^{N}
        \left(
        \sum_{m=1}^{N-j}
            Y^{(j+m)}_{s}
            \left(
                X^{m}_{s,s'}
            \right)
            +
            R^{(j)}_{s,s'}
        \right)
        \left(
            X^{j+1}_{s',t}
        \right)
    -
    \sum_{l=0}^{N-1}
        \sum_{j=l+1}^{N}
            Y_s^{(j)}
            \left(
                X_{s,s'}^{j-l}
            \right)
            \left(
                X_{s',t}^{l+1}
            \right)
    \\[0.5em]
    &=
    \sum_{j=0}^{N}
        \left(
        \sum_{m=1}^{N-j}
            Y^{(j+m)}_{s}
            \left(
                X^{m}_{s,s'}
            \right)
            +
            R^{(j)}_{s,s'}
        \right)
        \left(
            X^{j+1}_{s',t}
        \right)
    -
    \sum_{l=0}^{N-1}
        \sum_{m=1}^{N-l}
            Y_s^{(l+m)}
            \Big(
                X_{s,s'}^{m}
            \Big)
            \Big(
                X_{s',t}^{l+1}
            \Big)
    \\[0.5em]
    &=
    \sum_{j=0}^N
        R^{(j)}_{s,s'}
        \left(
            X^{j+1}_{s',t}
        \right),
\end{align*}
completing the proof.
\end{proof}

Define the space $C^{p, \beta}_{\omega} ([0,T]; E)$ to contain all $\Xi : [0,T]^2 \rightarrow E$ such that $\| \Xi \|_{p} + \| \delta \Xi \|_{\beta} < \infty$. We think of these $\Xi$ as describing a family of local approximations to some additive function; the map $\delta \Xi$ measures the failure of $\Xi$ to be additive itself. A key result is the sewing lemma, which is used to establish the existence, uniqueness and the regularity of this underlying additive function under suitable conditions on $\Xi$ and $\delta \Xi$. A central application of it is to provide the definition of and fundamental bounds on the rough integral of a controlled paths.

\begin{lemma}[Sewing lemma]\label{sewing lemma}
Let $p, \beta$ be such that $0 < \beta < 1 \leq p$. There exists a unique map $\Ical: C^{p, \beta}_{\omega} ([0,T]; E) \rightarrow C^{p}_{\omega}([0,T]; E)$ such that $(\Ical \Xi)_{0} = 0$ and
\begin{equation}
    \label{sewing map bound}
    \left|
        (\Ical \Xi)_{s,t}
        -
        \Xi_{s,t}
    \right|
    \leq
    \zeta(1/\beta)
    \;
    \omega(s,t)^{1/\beta}
    \|
        \delta \Xi
    \|_{\beta, \omega}
    ,
\end{equation}
where $\zeta$ is the Riemann-Zeta function and $(\Ical \Xi)_{s,t} = (\Ical \Xi)_t - (\Ical \Xi)_s$. Moreover $\Ical \Xi$ is such that
\[
    \Ical \Xi_{s,t}
    =
    \lim_{|\Dcal| \rightarrow 0}
        \sum_{m = 1}^{m_0}
            \Xi_{s_{m-1}, s_m}
\]
for partitions $\Dcal = \{s_0 < \dotsc < s_{m_0} \}$ over $[s,t]$.
\end{lemma}

The existence of the sewing map $\Ical$ was first presented in \cite{FdLP06} and is the core component of constructing rough integrals, whereby we take a sufficiently good local approximation and sew these together with the map $\Ical$. The proof of the lemma in the form given here follows a Young type argument as presented in Lemma 4.2 of \cite{FH14} with some modifications. The proof here also serves as a simpler example of the same strategy used to prove the existence of the joint rough integral in Theorem \ref{integral existence theorem}.

\begin{proof}
We begin this by establishing Young's maximal inequality in the one-parameter case. Let $\Dcal = \{s=s_0 < \dotsc < s_{m_0}=t\} \subset [s,t]$ and write
\[
    \sum_{\Dcal} \Xi
    =
    \sum_{m = 1}^{m_0}
        \Xi_{s_{m-1}, s_m}
    .
\]
The idea here is to selectively remove points from $\Dcal$ until we are left with the trivial partition $\{s,t\}$. Choose $m^{*}$ such that it minimises
\[
    \Big|
        \sum_{\Dcal}
            \Xi
        -
        \sum_{\Dcal\setminus \{s_m\}}
            \Xi
        \;
    \Big|
    =
    \left|
        \delta\Xi_{s_{m-1}, s_{m+1}}
    \right|
    .
\]
Then we have
\begin{align*}
    \Big|
        \sum_{\Dcal}
            \Xi
        -
        \sum_{\Dcal\setminus \{s_{m^*}\}}
            \Xi
        \;
    \Big|^{\beta}
    &\leq
    \frac{1}{m_0 - 1}
    \Big|
        \sum_{\Dcal}
            \Xi
        -
        \sum_{\Dcal\setminus \{s_m\}}
            \Xi
        \;
    \Big|^{\beta}
    \leq
    \frac{1}{m_0 - 1}
    \omega(s,t)
    \left\|
        \delta\Xi
    \right\|_{\beta, \omega}^{\beta}
    .
\end{align*}
Repeatedly choosing points $m^*$ to remove then gives us
\begin{align*}
    \Big|
        \sum_{\Dcal}
            \Xi
        -
        \Xi_{s,t}
        \;
    \Big|
    &=
    \sum_{l=1}^{m_0 - 1}
        \left(
            \frac{1}{m_0 - l}
            \omega(s,t)
            \left\|
                \delta\Xi
            \right\|_{\beta, \omega}^{\beta}
        \right)^{1/\beta}
    \leq
    \zeta
    \left(
        1/\beta
    \right)
    \omega(s,t)^{1/\beta}
    \left\|
        \delta\Xi
    \right\|_{\beta, \omega}
    .
\end{align*}
With this maximal bound on the sums over partitions, we now look to prove existence and uniqueness of the limit of these sums as the partition mesh size goes to zero. By definition the control $\omega$ is continuous and thus uniformly continuous on $[0, T]^2$. So for any $\epsilon > 0$ there exists $\delta_{\epsilon} > 0$ such that if $|s - t| < \epsilon$ then $|\omega(s,t)| < \delta_{\epsilon}$ and moreover $\delta_{\epsilon} \rightarrow 0$ as $\epsilon \rightarrow 0$.

Suppose that $\Dcal$ and $\Dcal'$ are two partitions of $[s,t]$, we wish to compare the sums of $\Xi$ over these two partitions as the mesh size decreases. Without loss of generality we will assume that $\Dcal \subset \Dcal'$, since otherwise we just compare the individual partitions with $\Dcal \cup \Dcal'$ for the same result. Writing $\Dcal = \{ s_0 < \dotsc < s_{m_0}\}$, we now view $\Dcal'$ as the union of partitions $\Dcal_{m}' = \Dcal' \cap [s_{m-1}, s_{m}] $ of $[s_{m-1}, s_{m}]$. Applying the earlier maximal bound:
\begin{align*}
    \Big|
        \sum_{\Dcal'}
            \Xi
        -
        \sum_{\Dcal}
            \Xi
        \;
    \Big|
    &\leq
    \sum_{m=1}^{m_0}
        \Big|
            \sum_{\Dcal_{m}'}
                \Xi
            -
            \Xi_{s_{m-1}, s_{m}}
            \;
        \Big|
    \\[0.5em]
    &\leq
    \sum_{m=1}^{m_0}
        \zeta
        \left(
            1/\beta
        \right)
        \omega(s_{m-1}, s_{m})^{1/\beta}
        \left\|
            \delta\Xi
        \right\|_{\beta, \omega}
    \\[0.5em]
    &\leq
    \zeta
    \left(
        1/\beta
    \right)
    \omega(s,t)
    \left\|
        \delta\Xi
    \right\|_{\beta, \omega}
    \delta_{\epsilon}^{1/\beta - 1}
    ,
\end{align*}
which goes to zero as $\epsilon$ goes to zero. Thus we have existence and uniqueness of the limit
\[
    (\Ical \Xi)_{s,t}
    =
    \lim_{|\Dcal| \rightarrow 0}
        \sum_{\Dcal \subset [s,t]}
            \Xi
    ,
\]
where $\Ical \Xi$ satisfies the bound \eqref{sewing map bound} and is such that $(\Ical \Xi)_0 = 0$.
\end{proof}

An elementary application of the sewing lemma on $\Xi^Y$ as defined in \eqref{local approximation} then yields existence of the rough integral as the limit of enhanced Riemann sums.

\begin{theorem}
Suppose that $\left(Y = Y^{(0)}, \dotsc, Y^{(N)}\right) \in \Dscr^p_X ([0,T]; E)$ with remainders denoted by $\left(R^{(0)}, \dotsc, R^{(N)} \right)$. Then for $[s, t] \subset [0, T]$, the rough integral
\[
    \int_s^t
        Y_r
    \; dX_{r}
    :=
    \lim_{|\Dcal| \rightarrow 0}
        \sum_{m = 1}^{m_0}
            \left(
                \sum_{j=0}^N
                    Y^{(j)}_{s_{m-1}}
                    \left(
                        X^{j+1}_{s_{m-1}, s_{m}}
                    \right)
            \right)
\]
exists and is such that
\begin{equation}
    \left|
        \int_s^t
            Y_r
        \; dX_{r}
        -
        \sum_{j=0}^N
            Y^{(j)}_{s}
            \left(
                X^{j+1}_{s,t}
            \right)
    \right|
    \leq
    \zeta\left(
        \theta
    \right)
    \omega(s,t)^{\theta}
    \sum_{j=0}^N
        \left\|
            X^{j+1}
        \right\|_{\frac{p}{j+1}}
        \left\|
            R^{(j)}
        \right\|_{\frac{p}{\pfloor - j}}
\end{equation}
where $\theta = \frac{\pfloor + 1}{p}$.
\end{theorem}

\end{document}